\documentclass[11pt]{article}%
\usepackage{natbib}
\usepackage{amsmath}
\usepackage{amssymb}
\usepackage{graphicx}
\usepackage{amsfonts}
\usepackage{color}%
\newtheorem{remark}{Remark}
\renewcommand{\cite}[2][]{\citep[#1]{#2}}
\begin{document}

\title{Morphing Ensemble Kalman Filters}
\author{Jonathan D. Beezley and Jan Mandel\\
Center for Computational
Mathematics, \\
University of Colorado at Denver and Health Sciences
Center, Denver, CO, \\  and \\
Mesoscale and Microscale Meteorology Division, \\ National Center for
Atmospheric Research, Boulder, CO }

\maketitle

\begin{abstract}
A new type of ensemble filter is proposed, which combines an ensemble
Kalman filter (EnKF) with the ideas of morphing and registration from
image processing. This results in filters suitable for nonlinear
problems whose solutions exhibit moving coherent features, such as
thin interfaces in wildfire modeling. The ensemble members are
represented as the composition of one common state with a spatial
transformation, called registration mapping, plus a~residual. A~fully
automatic registration method is used that requires only gridded
data, so the features in the model state do not need to be identified
by the user. The morphing EnKF operates on a transformed state
consisting of the registration mapping and the residual. Essentially,
the morphing EnKF uses intermediate states obtained by morphing
instead of linear combinations of the states.
\end{abstract}%

\section{Introduction}

The research reported here has been motivated by data assimilation
into wildfire models \cite{Mandel-2006-WMD}. Wildfire modeling
presents a~challenge to data assimilation because of non-Gaussian
probability distributions centered around the burning and not burning
states, and because of movements of thin reaction fronts with sharp
interfaces. This work is a part of an effort to build a Dynamic Data
Driven Application System (DDDAS, \citet{Darema-2004-DDD}) for
wildfires \cite{Mandel-2007-DDD}. Data assimilation is one of the
building blocks of the DDDAS concept, which involves a symbiotic
network of computer simulations and sensors.

The standard Ensemble Kalman Filter (EnKF) approach \cite{Evensen-2003-EKF}
fails for highly nonlinear problems that have solutions with coherent
features, such as firelines, rain fronts, or vortices, because it is limited
to a linear update of state. This can be ameliorated to some degree by
penalization of nonphysical states \cite{Johns-2005-CEK}, localization of EnKF
\cite{Anderson-2003-LLS,Ott-2003-LEK}, and employing the location of the
feature as an observation function \cite{Chen-2007-AVP}, but the basic problem
remains: EnKF works only when the increment in the location of the feature is
small \cite{Chen-2007-AVP}.
%due to model errors or infrequent data.

EnKF analysis formulas are based on the assumption that all
probability distribution involved are Gaussian, so even if an
ensemble can approximate non-Gaussian distribution, the closer the
system state distribution is to Gaussian the better. One mechanism
how non-Gaussian distributions arise in strongly nonlinear systems is
by evolution of coherent features. While the location and the size of
the feature may have an error distribution that is approximately
Gaussian, this is not necessarily the case for the value of the state
at a~given point. E.g., the typical probability density of the
temperature at a~point near the reaction region in a~wildfire model
is concentrated around the burning and the ambient temperature
(Fig.~\ref{fig:da_pdf}a). There is clearly a~need to adjust the
simulation state by distorting the simulation state in space rather
than employing an additive correction to the state. Therefore,
alternative error models that include the position of features were
considered in the literature \cite{Hoffman-1995-DRF,Davis-2006-OVP1}
and a~number of works emerged that achieve more efficient movement of
features by using a~spatial transformation as the field to which
additive corrections are made, such as a~transformation of the space
by a~global low order polynomial mapping to achieve alignment
\cite{Alexander-1998-UDW}, and two-step models to use alignment as
preprocessing to an additive correction
\cite{Lawson-2005-AEM,Ravela-2007-DAF}.

Moving and stretching one given image to become another given image is known
in image processing as registration \cite{Brown-1992-SIR}. Once the two images
are registered, one can easily create intermediate images, which is known as morphing.

The essence of the new method described here is to replace the linear
combinations of states in an ensemble filter by intermediate states obtained
by morphing. This method provides additive and position correction in a~single
step. For the analysis step (the Bayesian update), the state is transformed
into an extended state consisting of additive and position components. After
the analysis step, the state is converted back and advanced in time. The
purpose of this article is to demonstrate the potential of this approach.

The paper is organized as follows. In Section \ref{sec:assimilation}, we
briefly recall data assimilation by EnKF but we do not present the formulation
of the EnKF in detail and refer to the literature instead. In Section
\ref{sec:registration}, we describe image morphing and the automatic
registration algorithm used. Section \ref{sec:filter} contains the formulation
of the morphing EnKF. Numerical results on a wildfire model problem are
reported in Section \ref{sec:numerical}. Section \ref{sec:conclusion} contains
the conclusion and a~discussion of future extensions.

\section{Data Assimilation and Ensemble Filters}

\label{sec:assimilation}

The purpose of data assimilation is to estimate the system state using all
data available up to the current time. A discrete filter, considered here,
works by advancing in time a~probability distribution of the model state until
a~given \emph{analysis time}. At the analysis time, the probability
distribution, now called the \emph{prior} or the \emph{forecast}, is modified
by accounting for the data, which are considered to be observed at that time.
The new probability distribution, called the \emph{posterior} or the
\emph{analysis}, is given by the Bayes theorem,%
\begin{equation}
p_{a}(u)\varpropto p\left(  d|u\right)  p_{f}\left(  u\right)  ,
\label{eq:Bayes-dens}%
\end{equation}
where $\varpropto$ means proportional, $u$ is the model state, $p_{f}\left(
u\right)  $ is the forecast probability density, $p_{a}(u)$ is the analysis
probability density, $d$ is the data, and $p\left(  d|u\right)  $ is the data
likelihood. The data likelihood is the density of the probability that the
data value is $d$ assuming that the state is $u$. Assuming an additive data
error model, the data likelihood is found from the probability density
$p_{\varepsilon}$ of data error, which is assumed to be known (every
measurement must be accompanied by an error estimate), and from
an~\emph{observation function} $h$, by%
\begin{equation}
p\left(  d|u\right)  =p_{\varepsilon}\left(  d-h\left(  u\right)  \right)  .
\label{eq:data-error}%
\end{equation}
The value $h\left(  u\right)  $ of the observation function is what the
correct value of the data would be if the model state $u$ were exact. For a
tutorial on data assimilation, see \citet{Kalnay-AMD-2003}.

The well-known Kalman filter \cite{Kalman-1960-NAL} reduces the Bayesian
update (\ref{eq:Bayes-dens}) to linear algebra in the case when all
probability distributions are Gaussian. The Kalman filter must advance the
covariance of state, which is possible only when the model is linear, and it
is computationally very expensive. The EnKF
\cite{Evensen-1994-SDA,Houtekamer-1998-DAE} approximates the probability
distribution by the empirical measure $\frac{1}{N}\sum_{i=1}^{N}\delta_{u_{i}%
}$, where $u_{i}$ are members of an \emph{ensemble} of simulations states, and
$\delta_{u}$ denotes the Dirac delta measure concentrated at $u$.
 Each ensemble member is advanced in time between the Bayesian updates
independently. The EnKF approximates the mean and the covariance of the
forecast by the mean and the covariance of the ensemble, while still making
the assumption that all probability distributions are Gaussian. The EnKF works
by forming the analysis ensemble as \emph{linear combinations of the forecast
ensemble}, and the Bayesian update is implemented by linear algebra operating
on the matrix created from ensemble members as columns. This allows an
efficient implementation using high-performance matrix operations. We have
used the EnKF version from \citet{Burgers-1998-ASE}, but any other variant could
be used as well. For surveys of EnKF techniques, see \citet{Evensen-2007-DAE}
 and \citet{Tippett-2003-ESR}.

\section{Registration, Warping, and Morphing}

\label{sec:registration}

In this section, we build the tools from image processing that we are going to
use for the morphing EnKF later in Section \ref{sec:filter}. The
registration problem in image processing is to find a~spatial mapping that
turns two given images into each other \cite{Brown-1992-SIR}. Classical
registration requires the user to pick manually points in the two images that
should be mapped onto each other, and the registration mapping is interpolated
from that. Here we are interested in a~\emph{fully automatic registration
procedure that does not require any user input}. The specific feature or
objects, such as fire fronts or hurricane vortices, do not need to be
specified. The method takes as input only the pixel values in an image (i.e.,
gridded arrays). Of course, for the method to work well, the images being
registered should be sufficiently similar.

\subsection{Notation}

We will find it useful to use mappings as a convenient notation, so we review
few basics in an informal manner. The symbol \textquotedblleft$\mapsto
$\textquotedblright\ is read \textquotedblleft maps to\textquotedblright, and
we sometimes find it convenient to write $A:z\mapsto w$ instead of $A\left(
z\right)  =w$. The domain of a mapping $A$ is the set of its arguments $z$
such that $A\left(  z\right)  $ is defined. For two mappings $A$ and $B$,
their composition $A\circ B$ is defined by $A\circ B:z\mapsto A\left(
B\left(  z\right)  \right)  $; $I$ is the identity mapping defined by
$I:z\mapsto z$ (for all $z$ from its domain); and the inverse $A^{-1}$ of $A$
is a mapping such that $A^{-1}\circ A=I$ and $A\circ A^{-1}=I$. A real
function on a domain $D$ is also a mapping, namely one that maps $D$ into the
reals, $\mathbb{R}$. Given two mappings $A$ and $B$ from the same domain into
the same linear space (say, $\mathbb{R}$ or $\mathbb{R}^{k}$) and two scalars
$a$ and $b$, the linear combination of the mappings is defined by taking the
linear combination of their values,%
\[
aA+bB:z\mapsto aA\left(  z\right)  +bB\left(  z\right)  .
\]

\subsection{Image Registration}

Consider, for example, two grayscale images with intensities $u$ and $v$,
given as functions on some domain $D$ (such as a rectangle in the plane,
$\mathbb{R}^{2}$). For simplicity, assume that both $u$ and $v$ are equal to
some constant at and near the boundary of the domain $D$. In our application,
$u$ and $v$ will be temperature fields from two states of our
wildfire model, the fire will be inside the domain, and the temperature near
the boundary of the domain $D$ will be equal to the ambient temperature
(assumed to be the same everywhere). In image processing, $u$ and $v$ can be the
darkness levels of two photographs of objects with the same solid background.
The functions $u$ and $v$ will be also referred to as images in the following description.

The registration then becomes the problem to find two functions $T_{x}\left(
x,y\right)  $ and $T_{y}\left(  x,y\right)  $ such that the transformation of
the argument of $u$ by%
\begin{equation}
\left(  x,y\right)  \mapsto\left(  x+T_{x}\left(  x,y\right)  ,y+T_{y}\left(
x,y\right)  \right)  \label{eq:map}%
\end{equation}
transforms $u$ into a function approximately equal to $v$ on the domain $D$,
\begin{equation}
v\left(  x,y\right)  \approx u\left(  x+T_{x}\left(  x,y\right)
,y+T_{y}\left(  x,y\right)  \right)  , \label{eq:registered}%
\end{equation}
for all $\left(  x,y\right)  \in D$.

Define the mapping $T$ from $D$ to $\mathbb{R}^{2}$ by
\begin{equation}
T:\left(  x,y\right)  \mapsto\left(  T_{x}\left(  x,y\right)  ,T_{y}\left(
x,y\right)  \right)  . \label{eq:def-T}%
\end{equation}
Then (\ref{eq:registered}) can be written compactly as%
\[
v\approx u\circ\left(  I+T\right)  \text{ on }D,
\]
or%
\begin{equation}
v-u\circ\left(  I+T\right)  \approx0\text{ on }D. \label{eq:registered-T}%
\end{equation}
The mapping $I+T$ will be called the \emph{registration mapping}, and
the mapping $T$ will be called \emph{warping}. The reason
for writing the registration mapping as $I+T$ is that the zero warping $T=0$
is the neutral element of the operation of addition, and so linear
combinations of warpings have a~meaningful interpretation as blends of the
warpings. This will be important in the development of the morphing EnKF.

To avoid unnecessarily large and complicated warping, the warping $T$ should
be as also close to zero and as smooth as possible,%
\begin{equation}
T\approx0,\quad\nabla T\approx0, \label{eq:small-T}%
\end{equation}
where $\nabla T$ denotes the matrix of the first derivatives (the Jacobian matrix)
of $T$,%
\[
\nabla T=%
\begin{pmatrix}
\frac{\partial T_{x}}{\partial x} & \frac{\partial T_{x}}{\partial y}\\
\frac{\partial T_{y}}{\partial x} & \frac{\partial T_{y}}{\partial y}%
\end{pmatrix}
.
\]

In addition, we require that the registration mapping $I+T$ is one-to-one, so
the inverse $\left(  I+T\right)  ^{-1}$ exists. However, we do not require
that the values of $I+T$ are always in $D$ or the inverse $\left(  I+T\right)
^{-1}$ is defined on all of $D$, so $u\circ\left(  I+T\right)  $ and
$u\circ\left(  I+T\right)  ^{-1}$may not be defined on all of $D$. Therefore,
we consider all functions $u$, $v$, $u\circ\left(  I+T\right)  ,$
$u\circ\left(  I+T\right)  ^{-1}$, etc., extended on the whole $\mathbb{R}%
^{2}$ by the constant value of $u$ and $v$ on the boundary of $D$.

\subsection{Morphing}

Once the registration mapping $I+T$ is found, one can construct intermediate
functions $u_{\lambda}$ between $u_{0}$ and $u_{1}$ by \emph{morphing}
(Fig.~\ref{fig:morph}),
\begin{equation}
u_{\lambda}=\left(  u+\lambda r\right)  \circ\left(  I+\lambda T\right)
,\quad0\leq\lambda\leq1, \label{eq:ulambda}%
\end{equation}
where
\begin{equation}
r=v\circ\left(  I+T\right)  ^{-1}-u \label{eq:inverse_res}%
\end{equation}
will be called the \emph{registration residual}; it is easy to see that $r$ is
linked to the approximation in (\ref{eq:registered-T}) by%
\[
r=\left(  v-u\circ\left(  I+T\right)  \right)  \circ\left(  I+T\right)
^{-1},
\]
thus (\ref{eq:registered-T}) also implies that $r\approx0$.

The original functions $u$ and $v$ are recovered by choosing in
(\ref{eq:ulambda}) $\lambda=0$ and $\lambda=1$, respectively,
\begin{align}
u_{0}  &  =u\circ I=u,\label{eq:back}\\
u_{1}  &  =\left(  u+r\right)  \circ\left(  I+T\right) \label{eq:back-1}\\
&  =\left(  u+v\circ\left(  I+T\right)  ^{-1}-u\right)  \circ\left(
I+T\right) \nonumber\\
&  =v\circ\left(  I+T\right)  ^{-1}\circ\left(  I+T\right)  =v.\nonumber
\end{align}

\begin{remark}
In the registration and morphing literature, the residual is often neglected.
Then the morphed function is given simply by the transformation of the
argument, $u_{\lambda}=u\circ\left(  I+\lambda T\right)  $. The simplest way
how to account for the residual is to add a correction term to $u_{\lambda}$.
This gives the morphing formula
\begin{equation}
u_{\lambda}=u\circ\left(  I+\lambda T\right)  +\lambda\left(  v-u\circ\left(
I+T\right)  \right)  ,\quad0\leq\lambda\leq1, \label{eq:simple-morph}%
\end{equation}
which is much easier to use because does not require the inverse $\left(
I+T\right)  ^{-1}$ like (\ref{eq:ulambda}). The formula (\ref{eq:simple-morph}%
) also recovers $u=u_{0}$ and $v=u_{1}$, but, in our computations, we have found
it unsatisfactory for tracking features and therefore we do not use it. The
reason is that when the residual is not negligibly small, the intermediate
functions $u_{\lambda}$ will have a spurious feature in the fixed location
where the residual is large. On the other hand, the more expensive improved
morphing formula (\ref{eq:ulambda}) moves the contribution to the change in
amplitude along with the change of the position.
\end{remark}

\subsection{Grids}

\label{sec:grids}An array of values associated with a rectangular grid is
called a \emph{gridded array}. The functions $u_{\lambda}$ are represented by
a gridded array on a \emph{pixel grid}, while the mapping $T$ is represented
by two gridded arrays on a coarser \emph{morphing grid}. In our application,
the domain $D$ is a rectangle, discretized by a uniform $n_{x}\times n_{y}$
pixel grid with $n=n_{x}n_{y}$ nodes, and, for $i=1,\ldots,M$, by a
$(2^{i}+1)\times(2^{i}+1)$ uniform grid $D_{i}$ with nodes denoted by $\left(
x_{j},y_{k}\right)$ and total
$m_{i}=\left(  2^{i}+1\right)  ^{2}$ nodes. The morphing grid is the finest
grid $D_{M}$, with $m=m_{M}$ nodes. Denote by $T_{i}$ the gridded array $T$
restricted to the grid $D_{i}$, and by $I_{i-1}^{i}$ the bilinear
interpolation operator from grid $D_{i-1}$ to grid $D_{i}$. All grids contain
nodes on the boundary of $D$. It is assumed that $n\gg m$.

The values of the functions and of the mappings away from their respective
grid points are evaluated by bilinear interpolation without mentioning the
interpolation explicitly. So, composed functions like $u\circ\left(
I+T\right)  $ are calculated in a straightforward manner: for an arbitrary
$\left(  x,y\right)  \in D$, first $(I+T)\left(  x,y\right)  =\left(
x^{\prime},y^{\prime}\right)  $ is computed by bilinear interpolation on the
morphing grid and then $u\left(  x^{\prime},y^{\prime}\right)  $ is evaluated
by bilinear interpolation on the pixel grid. The calculation of $(I+T)^{-1}%
\left(  x^{\prime},y^{\prime}\right)  $ is done by inverse interpolation as
follows. For uniformly spaced nodes $\left(  x_{j},y_{k}\right)  $ on the
morphing grid, the images $\left(  x_{j}^{\prime},y_{k}^{\prime}\right)
=(I+T)\left(  x_{j},y_{k}\right)  $ form a nonuniform grid, and the value of
$(I+T)^{-1}(x^{\prime},y^{\prime})$ is approximated by interpolating the
original coordinates $x$ and $y$ as functions on the nonuniform grid. This can
be accomplished efficiently, e.g., by using the MATLAB function
\texttt{griddata}, which works exactly like interpolation, but allows for
nonuniformly spaced data. We then have $\left(  I+T\right)  ^{-1}\circ\left(
I+T\right)  \approx I$ and $u_{1}\approx v$ as in (\ref{eq:back-1}) only up to
the error of interpolation from the morphing grid.

\subsection{An Automatic Registration Procedure}

The formulation of registration as (\ref{eq:registered-T}) --
(\ref{eq:small-T}) naturally leads to a construction of the mapping $T$ by
optimization. So, suppose we are given $u$ and $v$ on the pixel grid and wish
to find a warping $T$ that is an approximate solution of%
\begin{equation}
J(T,u,v)=\Vert v-u\circ\left(  I+T\right)  \Vert+C_{1}\Vert T\Vert+C_{2}%
\Vert\bigtriangledown T\Vert\rightarrow\min_{T}, \label{eq:registration}%
\end{equation}
where the norms are chosen as%
\begin{align}
&  \Vert v-u\circ\left(  I+T\right)  \Vert=\int_{D}\left\vert v-u\circ\left(
I+T\right)  \right\vert dxdy,\label{eq:norm-r}\\
&  \Vert T\Vert=\int_{D}\left\vert T_{x}\right\vert +\left\vert T_{y}%
\right\vert dxdy,\label{eq:norm-T}\\
&  \Vert\bigtriangledown T\Vert=\int_{D}\left\vert \frac{\partial T_{x}%
}{\partial x}\right\vert +\left\vert \frac{\partial T_{x}}{\partial
y}\right\vert +\left\vert \frac{\partial T_{y}}{\partial x}\right\vert
+\left\vert \frac{\partial T_{y}}{\partial y}\right\vert dxdy.
\label{eq:norm-T-grad}%
\end{align}

The optimization formulation tries to balance the conflicting objectives of
good approximation by the registered image, and as small and smooth warping as
possible. The objective function $J(T,u,v)$ is in general
not a convex function of $T$, and so there are many local minima. For example,
a local minimum of may occur when some small part of $u$ and $v$ matches,
while the overall match is still not very good.

To solve the minimization problem (\ref{eq:registration}), we have used the
algorithm from \citet{Gao-1998-WMA} with several modifications. We have also
filled in some details that were not provided. The description of the
resulting algorithm is the subject of this section. It is presented for
completeness only; any other automatic registration procedure from the
literature could be used as well. The algorithm guarantees by construction
that $I+T$ is invertible.

For speed and to decrease the chance that the minimization gets stuck in
a~local minimum, the method proceeds by building $T$ on a~nested hierarchy of
meshes $D_{i}$, starting with $T_{1}$ on the coarsest mesh $D_{1}$ and ending
with the mapping $T=T_{M}$ built on the morphing grid $D_{M}$. In the
computation, $\|v-u\circ\left(  I+T_{i}\right)\|  $ in (\ref{eq:norm-r}) is
integrated numerically on the pixel grid, while $\Vert T\Vert$and
$\Vert\bigtriangledown T\Vert$ in (\ref{eq:norm-T}) and (\ref{eq:norm-T-grad})
are integrated on one of the grids $D_{i}$, with the derivatives approximated
by finite differences. The integrals are approximated by the scaled sums of
the values of the integrands at grid nodes of the respective grids. Denote
this version of the objective function on the grid $D_{i}$ by $J_{i}%
(T_{i},u,v)$. Assume that an initial guess $\widetilde{T}$ of $T$ on the
morphing grid is known. If none is given, use $\widetilde{T}=0$. The notation
$\widetilde{T}_{i}$ means the gridded array $\widetilde{T}$ restricted to the
grid $D_{i}$, just like $T_{i}$.

On each mesh $D_{i}$, the method proceeds as follows. In order not to overload
the notation with many iteration indices, the values of $T$, and thus also
$T_{i}$, change during the computation just like a variable in a computer program.

\begin{enumerate}
\item Smooth $u$ and $v$ by convolution to get $\tilde{u}_{i}$ and $\tilde
{v}_{i}$ on the pixel grid.

\item If $i=1$, initialize $T_{1}=\widetilde{T}_{1}$. Otherwise interpolate
$T_{i}$ from the coarse grid $T_{i-1}$ as a correction to $\widetilde{T}$.

\item Minimize the objective function $J_{i}(T_{i},\tilde{u}_{i},\tilde{v}%
_{i})$ by adjusting the value of $T$ at one node of $D_{i}$ at a time. Stop
when a maximum number of sweeps through all nodes has been reached, or a
stopping test based on the decrease of the objective function or residual size
has been met.
\end{enumerate}

We now describe each step in more detail.

\begin{enumerate}
\item \emph{Smoothing. }The purpose of registration on a coarse mesh first is to
capture coarse similarities between the images $u$ and $v$. In order to force
coarse grids to capture coarse features only and to disregard fine features,
on each grid, we first smooth the images by convolution with a Gaussian
kernel. This allows to track large scale perturbations on coarse
grids even for a thin feature such as a fireline, while maintaining small
scale accuracy on fine grids. For the use with the grid $D_{i}$, we create
the smoothed image $\tilde{u}_{i}$ with resolution on the scale of $D_{i}$, by%
\begin{equation}
\tilde{u}_{i}\left(  x_{j},y_{k}\right)  =\sum_{j^{\prime}=-n_{x}+1}^{2n_{x}%
}\sum_{k^{\prime}=-n_{y}+1}^{2n_{y}}\hat{\phi}_{j}\left(  x_{j^{\prime}%
}\right)  u\left(  x_{j^{\prime}},y_{k^{\prime}}\right)  \hat{\psi}_{k}\left(
y_{k^{\prime}}\right)  , \label{eq:smoothing-conv}%
\end{equation}
where
\begin{align*}
\hat{\phi}_{j}\left(  x\right)   &  =c_{j}\exp\left(  -\frac{\left(
x_{j}-x\right)  ^{2}}{\alpha_{i}}\right)  ,\\
\hat{\psi}_{k}\left(  y\right)   &  =d_{k}\exp\left(  -\frac{\left(
y_{j}-y\right)  ^{2}}{\alpha_{i}}\right)  ,
\end{align*}
the constant $\alpha_{i}=0.25/(2^{i}+1)$ is tuned so that there is more
smoothing on the coarser grids, the values $u\left(  y_{j^{\prime}%
},y_{k^{\prime}}\right)  $ outside of $D$ are replaced by the constant
boundary value, and the normalization constants $c_{j}$ and $d_{k}$ are
determined so that%
\[
\sum_{j^{\prime}=-n_{x}+1}^{2n_{x}}\hat{\phi}_{j}^{2}\left(  x_{j^{\prime}%
}\right)  =\sum_{k^{\prime}=-n_{y}+1}^{2n_{y}}\hat{\psi}_{k}^{2}\left(
y_{k^{\prime}}\right)  =1.
\]
We also compute $\tilde{v}_{i}$ from $v$ in the same way.

\item \emph{Initialization.} Consider the grid $D_{i}$, $i>1$, with the nodes
$\left(  x_{j},y_{k}\right)  $. The values of $I+T_{i}$ are already known at
the nodes on the coarse grid $D_{i-1}.$ The correction that the optimization on
grid $D_{i-1}$ applied to the initial guess is $T_{i-1}-\tilde{T}_{i-1}$. To
apply this correction to $\tilde{T}$ everywhere on $D_{i}$, we interpolate the
correction from the grid $D_{i-1}$ to $D_{i}$ by the bilinear interpolation
$I_{i-1}^{i}$ to get the initial guess on the grid $D_{i}$, i.e.,%
\[
T_{i}=\tilde{T}_{i}+I_{i-1}^{i}\left(  T_{i-1}-\tilde{T}_{i-1}\right)  .
\]

\item \emph{Optimization.} The value of $A^{\prime}=\left(  I+T\right)
\left(  x_{j},y_{k}\right)  $ is optimized by first evaluating the local
objective function $f\left(  A^{\prime}\right)  =J_{_{i}}(T_{i},\tilde{u}%
_{i},\tilde{v}_{i})$ at the nodes of a local grid inside the mapped local
rectangle $\left(  I+T_{i}\right)  \left(  \left[  x_{j-1},x_{j+1}\right]
\times\left[  y_{k-1},y_{k+1}\right]  \right)  $. The location of $A^{\prime}$
is then refined by several iterations of nonlinear optimization, starting from
the local grid node with the least value of $f\left(  A^{\prime}\right)  $
found. We have used coordinate descent alternating in the $x$ and $y$
direction by calling MATLAB function \texttt{minfnb} for 1D constrained
optimization. Cf., Fig. \ref{fig:search}. For nodes $\left(  x_{j}%
,y_{k}\right)  $ on the boundary of the domain $D$, the location of
$A^{\prime}$ is constrained within $D$, but it is allowed to move inside the
domain $D$.
\end{enumerate}

The differences between the method described here and the method by
\citet{Gao-1998-WMA} are: the refinement of node positions by
nonlinear optimization; the gradient term in the objective function;
smoothing of the images before registration on coarse levels; the use
of an initial guess $\tilde{T}$; and allowing the nodes on the
boundary to move inside of the domain.

\subsection{Computational Complexity}

The operation (\ref{eq:smoothing-conv}) is the multiplication of three dense
matrices. Assuming bounded aspect ratio of the image, $n_{x}\approx
\operatorname*{const}n_{y}\approx\operatorname*{const}n^{0.5}$, the
computation (\ref{eq:smoothing-conv}) takes $O\left(  n^{1.5}\right)  $
operations. Using FFT to replace the convolution of functions by
multiplication of their Fourier coefficients, it can be implemented in
$O\left(  n\log n\right)  $ operations.

The cost of evaluating the entire objective function is linear in the number
of pixels in the image, which is not practical. Fortunately, computing the
entire objective function is unnecessary: changing $T\left(  x_{i}%
,y_{k}\right)  $ on the grid $D_{i}$ can only influence the terms in the
objective function associated with the region $D_{jk}=\left[  x_{j-1}%
,y_{k-1}\right]  \times\left[  x_{j+1},y_{k+1}\right]  $, which requires only
$O\left(n/\left(  2^{i}+1\right)^2  \right)$ operations. Since there are
$\left(  2^{i}+1\right)^2 $ nodes to
optimize, the cost of one optimization sweep is $O\left(  n\right)  .$

Recall that $m=(2^{M}+1)^{2}$ is the number of the morphing grid points. Since
the optimization on each grid cost $O\left(  n\right)  $ operations, smoothing
on each grid costs $O(n\log n)$ operations, and there is $M=O\left(  \log
m\right)  $ grids, the total complexity of the registration algorithm is
$O\left(  n\log m\log n\right)  $. Thus, the method is suitable for a large
number of pixels as well as a large number of nodes on the morphing grid.

\section{Morphing Ensemble Filter}

\label{sec:filter}

The state of the model in general consists of several gridded arrays,
$U=\left(  w,z,\ldots\right)  $. For simplicity, suppose that all arrays are
defined over the same grid and that the registration is applied only to the
first array, $w$; this will be the case in the model application in Section
\ref{sec:numerical}. The general case will be discussed in Section
\ref{sec:conclusion}.

Let $\left\{  U_{k}\right\}  =\left\{  U_{1},\ldots,U_{N}\right\}  $ be an
ensemble of states, with the ensemble member $U_{k}$ consisting of the gridded
arrays
\[
U_{k}=\left(  w_{k},z_{k},\ldots\right)  .
\]
The subscript $_{k}$ in this section means the number of the ensemble member,
and it is not associated with the hierarchy of grids as in Section
\ref{sec:registration}. The concept of the hierarchy of grids is relevant only
to the internal working of the particular automatic registration algorithm
described in Section \ref{sec:registration}; here we just use the result of
the registration, which is a mapping defined by gridded arrays on the morphing
grid $D_{M}$.

Given one fixed state $U_{0}=\left(  w_{0},z_{0},\ldots\right)  $, the
automatic registration (\ref{eq:registration}) of the first array $w$ defines
the \emph{registration representations} $\left[  R_{k},T_{k}\right]  $ of the
ensemble members as morphs of $U_{0}$, with the registration residual%
\begin{align*}
R_{k}  &  =\left(  r_{w_{k}},r_{z_{k}},\ldots\right)  ,\\
r_{w_{k}}  &  =w_{k}\circ\left(  I+T_{k}\right)  ^{-1}-w_{0},\\
r_{z_{k}}  &  =z_{k}\circ\left(  I+T_{k}\right)  ^{-1}-z_{0},\\
&  \vdots
\end{align*}
and warpings $T_{k}$ determined as approximate solutions of independent
optimization problems based on the state array $w$,
\[
\left\Vert w_{k}-w_{0}\circ\left(  I+T_{k}\right)  \right\Vert +\left\Vert
T_{k}\right\Vert +\left\Vert \nabla T_{k}\right\Vert \rightarrow\min_{T_{k}}.
\]
The mapping $T_{k}$ from the previous analysis cycle is used as the initial
$\widetilde{T}_{k}$ in the automatic registration. In our tests, this all but
guarantees good registration and a speedup of one or more orders of magnitude
compared to starting from zero.

Instead of EnKF operating on the ensemble $\left\{  U_{k}\right\}  $ and making
linear combinations of its members, the morphing EnKF applies the EnKF
algorithm to the ensemble of registration representations $\left\{  \left[
R_{k},T_{k}\right]  \right\}  $, resulting in the analysis ensemble in registration
representation, $\left\{  \left[  R_{k}^{a},T_{k}^{a}\right]  \right\}  $,
with $R_{k}^{a}=\left(  r_{w_{k}}^{a},r_{z_{k}}^{a},\ldots\right)  $. The
analysis ensemble is then transformed back by (\ref{eq:back-1}), which here
becomes%
\begin{align}
w_{k}^{a}  &  =\left(  w_{k}+r_{w_{k}}^{a}\right)  \circ\left(  I+T_{k}%
^{a}\right)  ,\label{eq:back-many}\\
z_{k}^{a}  &  =\left(  z_{k}+r_{z_{k}}^{a}\right)  \circ\left(  I+T_{k}%
^{a}\right)  ,\label{eq:back-z}\\
&  \vdots\nonumber
\end{align}

\begin{remark}
Note that the registration representations $\left\{  \left[  R_{k}^{a}%
,T_{k}^{a}\right]  \right\}  $ of the analysis ensemble are linear
combinations of the registration representations $\left\{  \left[  R_{k}%
,T_{k}\right]  \right\}  $ of the forecast ensemble. Denote by $\lambda_{k}$
the coefficients of one such linear combination; then a~member of the analysis
ensemble has the form%
\begin{equation}
\left(  w_{0}+\sum_{k=1}^{N}\lambda_{k}r_{w_{k}}\right)  \circ\left(
I+\sum_{k=1}^{N}\lambda_{k}T_{k}\right)  , \label{eq:analysis}%
\end{equation}
and similarly for the other state arrays. This imposes certain constraints,
e.g., in general it may not be possible to write the zero state as
(\ref{eq:analysis}), and thus the potential for amplitude corrections might be
limited. This limitation does not seem to be important in the application of
interest here (wildfire), and its effect will be studied elsewhere.
\end{remark}

Given an observation function $h$, cf., (\ref{eq:data-error}), the transformed
observation function for EnKF on the registration representations can be
obtained directly by substituting from (\ref{eq:back-many}) into the
observation function,
\[
\tilde{h}\left(  \left[  R,T\right]  \right)  =h\bigl(\left(  w+r_{w}\right)
\circ\left(  I+T\right)  ,\left(  z+r_{z}\right)  \circ\left(  I+T\right)
,\ldots\bigr).
\]
However, constructing the observation function this way may not be the best.
Consider the case of one point observation, such as the temperature at some
location. Then the difference between the observed temperature and the value
of the observation function gives little indication which way should the
transformed state be adjusted. Suppose the temperature reading is high and the
ensemble members have high temperature only in some small location (fireline).
Then it is quite possible that the observation function (temperature at the
given location) evaluated on ensemble members will miss the fire in the
ensemble members completely. This is, however, a reflection of the inherent
difficulty of localizing small features from point observations.

For data that is given as gridded arrays (e.g, images, or a dense array of
measurements), there is a better way. Suppose the data $d$ is a measurement of
the first array in the state, $w$. Then, transforming the data $d$ into its
registration representation $\left[  r_{d},T_{d}\right]  $ just like the
registration of the state array $w$, the observation equation
(\ref{eq:data-error}) becomes the comparison between the registration
representations of the data $d$ and the state array $w$,
\begin{equation}
\tilde{h}\left(  \left[  R,T\right]  \right)  =\left[  r_{d},T_{d}\right]
\approx\left[  r_{w},T\right]  . \label{eq:morph-obs}%
\end{equation}
Data given on a part of the domain can be registered and used in the same way.
Note that no manual identification of the location of the feature either in
the data or in the ensemble members is needed.

\section{Numerical Results}

\label{sec:numerical}

We have applied the morphing EnKF to an ensemble from the wildland
fire model in \citet{Mandel-2006-WMD}. The simulation state consists
of the temperature and the fuel fraction remaining on
a~$500\,m\times500\,m$ domain with a~$2\,m$ uniform grid. The model has two
state arrays, the temperature $w$ and fuel supply $z$. An initial
fire solution $U_{0}$ was created by raising a~small square in the center of
the domain above the ignition temperature and applying a~small amount
of ambient wind until a fire front developed. The simulated data consisted of the whole
temperature field of another solution, started from $U_{0}$ and
morphed so that the fire was moved to a different location of the
domain compared with the average location of the ensemble members.
The observation equation (\ref{eq:morph-obs}) was used, with Gaussian
error in the registration residual of the temperature and in the registration
mapping. The standard deviation of the data error was $50\,K$ and $5\,m$,
respectively. This large discrepancy is used to show that the
morphing EnKF can be effective even when data is very different from
the ensemble mean. The image registration algorithm was applied with
a $17\times17$ morphing grid, i.e., $M=4$ refinement levels. We have
performed up to 5 optimization sweeps, stopping if the relative
improvement of the objective function for the last sweep was less
than $0.001$ or if the infinity norm of the residual $r_{w}$ fell
below $1$ $K$. The optimization parameters used for scaling the norms
in the objective function (\ref{eq:registration}) were $C_1=10000\,
mK^{-1}$ and $C_2=1000\,m^2K^{-1}$. For simplicity in the
computation, the fuel supply variables were not included in the data
assimilation.  Although the fuel supply was warped spatially as in
(\ref{eq:back-z}), the registration residual of the fuel supply, $r_{z_k}$, was taken to
be zero.

The $50$ member ensemble shown in Fig.~\ref{fig:da_morph}
was generated by morphing the initial state $U_{0}$ using smooth random
fields $r_{w_{k}}$ and $T_{k}$ of the form
\begin{equation}
\operatorname{const}\sum_{j=1}^{d}\sum_{\ell=1}^{d}\lambda_{j,\ell}d_{j,\ell}\sin j\pi x\sin
\ell\pi y
\label{eq:smooth-field}
\end{equation}
with $d_{j,\ell}\sim N(0,1)$ and $\lambda_{j,\ell}=\left(  1+\sqrt{j^{2}%
+\ell^{2}}\right)  ^{-2}$, for each ensemble member. The constant in
(\ref{eq:smooth-field}) was $50K$ for the residual $r_{w_k}$ and $5m$ for the warping
 $T_{k}$. Since it is not guaranteed that $(I+T)^{-1}$ exists for a smooth random $T$,
 we have tested if $I+T$ is one to one and generated another random $T$ if not.
The resulting $17\times17$ and $250\times250$ matrices are appended to form
$17^{2}+250^{2}$ element vectors representing an ensemble state $[R_{k}%
,T_{k}]$ for the EnKF. The same state $U_{0}$ was advanced in time along with the
ensemble. (Of course, other choices of $U_{0}$ are possible.)

The ensemble was advanced in $3$ minute analysis cycles. The new
registration representation $[r_{k},T_{k}]$ was then calculated using
the previous analysis values as an initial guess and incremented by
EnKF. The ensemble after the first analysis cycle is shown in
Fig.~\ref{fig:da_morph}. The results after five analysis cycles were
similar, indicating no filter divergence (Fig.~\ref{fig:da_morph5}).
Numerical results indicate that the error distribution of the
registration representation is much closer to Gaussian than the error
distribution of the temperature itself. This is demonstrated in
Fig.~\ref{fig:da_pdf}, where the estimated probability density
functions for the temperature, the registration residual of the temperature,
and the registration mapping for
Fig.~\ref{fig:da_morph5}
are computed at a single point in the domain using a Gaussian kernel with bandwidth
$0.3$ times the sample standard deviation. The point was chosen to be on the boundary of the
fire shape, where the non-Gaussianity may be expected to be the strongest.
In Fig.~\ref{fig:da_normtest5}, the Anderson-Darling test for normality
was applied to each point on the domain for the analysis step from Fig.~\ref{fig:da_morph5}.
The resulting $p$-values
were plotted on their corresponding locations in the domain with
darkness determined on a log scale with black as $10^{-8}$ (highly
non-Gaussian) and white $1$ (highly Gaussian). While the
Anderson-Darling test is intended to be a hypothesis test for
normality, it is used here to visualize on a continuous scale the
closeness to normality of the marginal probability
distribution any point in the domain. Again, strongest departure from normality of the
distribution is seen around the fire.

The implementation was a prototype done in Matlab. Therefore, we do
not report timings.

\section{Conclusion}

\label{sec:conclusion}

The numerical results show that the morphing EnKF is useful for a a~highly
nonlinear problem (a~model problem for wildfire simulation) with a coherent spatial feature of
the solution (propagating fireline). In previous work
\cite{Johns-2005-CEK,Mandel-2006-WMD}, we have used penalization of
nonphysical solutions, but the location of the fire in the data could not be
too far from the location in the ensemble, artificial perturbation had to be
added to retain the spread of the ensemble, and the penalty constant, the
amount of additional spread, and the data variance had to be finely tuned for
acceptable results. This new method does not appear to have the same
limitations. The registration works automatically from gridded data and no
objects need to be specified. The difference between the feature location in
the data and in the ensemble can be large and the data variance can be as
small as necessary, without causing filter divergence. One essential
limitation is that the registered images need to be sufficiently similar, and
the registration mapping should be sufficiently similar to the initial guess.
This will eventually impose a limitation on how long can the ensemble go
without an analysis step. However, compared to previous results for the same
problem \cite{Johns-2005-CEK,Mandel-2006-WMD}, the convergence of the present
filter is much better.

It was shown that the number of operations grows almost linearly with
the number of degrees of freedom, so the method is suitable for very
large problems. The method may be useful in wildfire modeling as well
as in data assimilation for other problems with strongly non-Gaussian
distributions and moving coherent features, such as rain fronts or
hurricane vortices. This paper only presents the basic method and
reports on results for a highly nonlinear, but still a fairly simple
model problem. Further enhancements, needed for practical
applications, such as general atmospheric science problems and
coupled wildfire and atmosphere models, will be studied elsewhere.
These enhancements should include the following.

In registration (Section \ref{sec:registration}), small residual can
be forced on a smaller domain than the whole domain to allow a global
shift and rotation of the image. This will be important in weather
applications, where there is no solid background; the ambient
temperature served as the background in the wildfire model tested
here. The registration method guarantees that the inverse $\left(
I+T\right)  ^{-1}$ exists, but the derivatives of the inverse could
be very large, resulting in a~loss of stability. Therefore, the
inverse should be involved in the objective function. For example,
penalty functions can be used to force in the local optimization step
the value of $T\left(  x_{j},y_{k}\right)  $ to stay well inside the
region where $\left( I+T\right)  ^{-1}$ exists, or the reciprocal of
the Jacobian of $I+T$ or the norm of the inverse $(I+T)^{-1}$,
multiplied by a constant, can be added as a penalty term to the
objective function directly. Also, the registration does not treat
the input images $u$ and $v$ in the same way; an algorithm that is
symmetric with respect to swapping the input images would take care
of the issue automatically. Finally, the registration mapping is
piecewise bilinear, and therefore the morphed state will have kinks
-- not very good for differential equations -- unless the morphing
grid is as fine as the pixel grid. In order to be able to use a
relatively coarse morphing grid (and thus cheap registration), one
could use smooth interpolation, such as $B$-splines, used in the
image registration context, e.g., by \citet{Arganda-2006-CER}.

The registration might be further improved by updating the
registration mapping more often, whether there are new data or not,
and thus assuring that the initial guess for the registration
algorithm is always sufficiently close. Also, the common state $U_0$
to register the ensemble members against has been evolved from one
initial condition regardless of the data. Over time, such $U_0$ could
diverge significantly from the ensemble (which tracks the data),
resulting in more strenuous registration. If this becomes an issue, a
better $U_0$ might be constructed from the analysis directly (perhaps
as the mean of registration representations of the analysis ensemble
members) and then evolved in time until the next analysis step.

When some ensemble members totally miss the feature (e.g., the fire),
the registration mapping does not matter much and all error will be
in the registration residual. This is not a problem, because those
ensemble members have low data likelihood, so they do not influence the
posterior pdf much. They do affect the ensemble mean and covariance, so the EnKF
analysis might change significantly if too many members miss the feature. Currently, the method was tested
for the case when there is a single significant feature (the fire),
which is essentially characterized by its position and strength.
Further research will be needed to deal with more general cases. The
method will need to be generalized to use more state arrays for the
registration at the same time, work with nested grids, and perhaps
use different registration mappings applied to different fields.
E.g., in a coupled atmosphere-fire model, the state of the atmosphere
and the state of the fire might require different position
adjustments. 3D registration might be needed for atmospheric
problems, especially those with strong buoyancy as over a wildfire.

\bibliographystyle{ametsoc}
\bibliography{../../../bibliography/dddas-jm}

\onecolumn

%all figures should be at the end

\newcommand{\captionexampleoned}{Example of morphing procedure
(\ref{eq:ulambda}) in one dimension. In this example, all variables
are nondimensional. The functions $u_0(x)$ and $u_1(x)$ are given on
the interval $[0,1]$. The intermediate functions $u_\lambda$ are
computed by the morphing formula (\ref{eq:ulambda}), which is seen to
interpolate the difference in position as well as the difference in
magnitude. The registration mapping $T$ is piecewise linear on the
morphing grid in the interval $[0,1]$ with spacing equal to $0.2$.
The thin dash-dot lines in the horizontal plane connect the nodes $x$
of the morphing grid with the values of the registration mapping
$(I+T)(x)$.}
\begin{figure}
\begin{center}%
\includegraphics[width=5in]{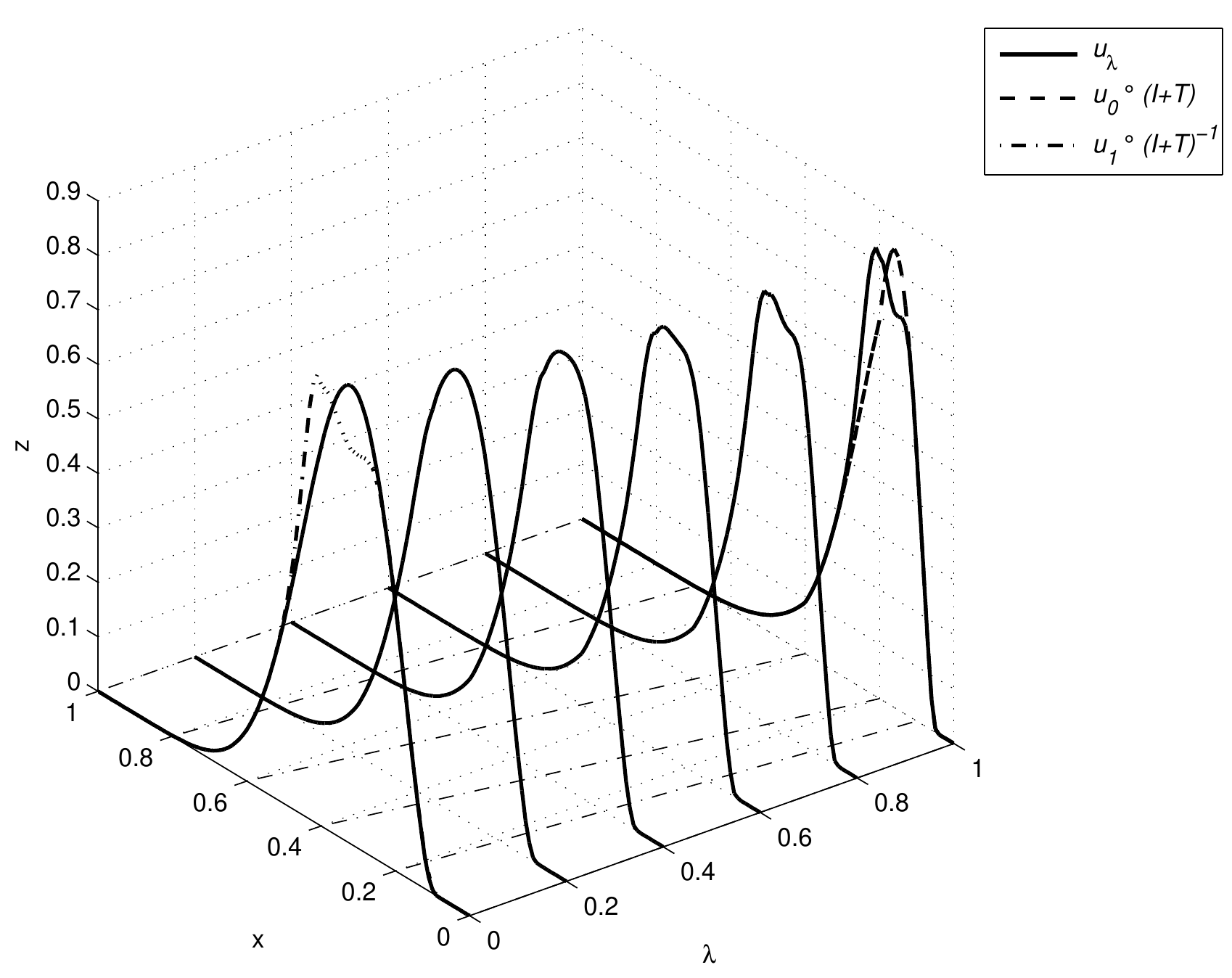}
\end{center}
\caption{\captionexampleoned}%
\label{fig:example1d}%
\end{figure}

\newcommand{\captionmorph}{Morphing of two solutions of a~reaction-diffusion equation system
used in a~wildfire simulation. The states with $\lambda=0$ and
$\lambda=1$ are given. The intermediate states are created
automatically. The horizontal plane is the earth surface. The
vertical axis and the color map are the temperature. The morphing
algorithm combines the values as well as the positions.}
\begin{figure}
\begin{center}%
\begin{tabular}
[c]{ccccc}%
\hspace*{-1em}\includegraphics[width=1.4in]{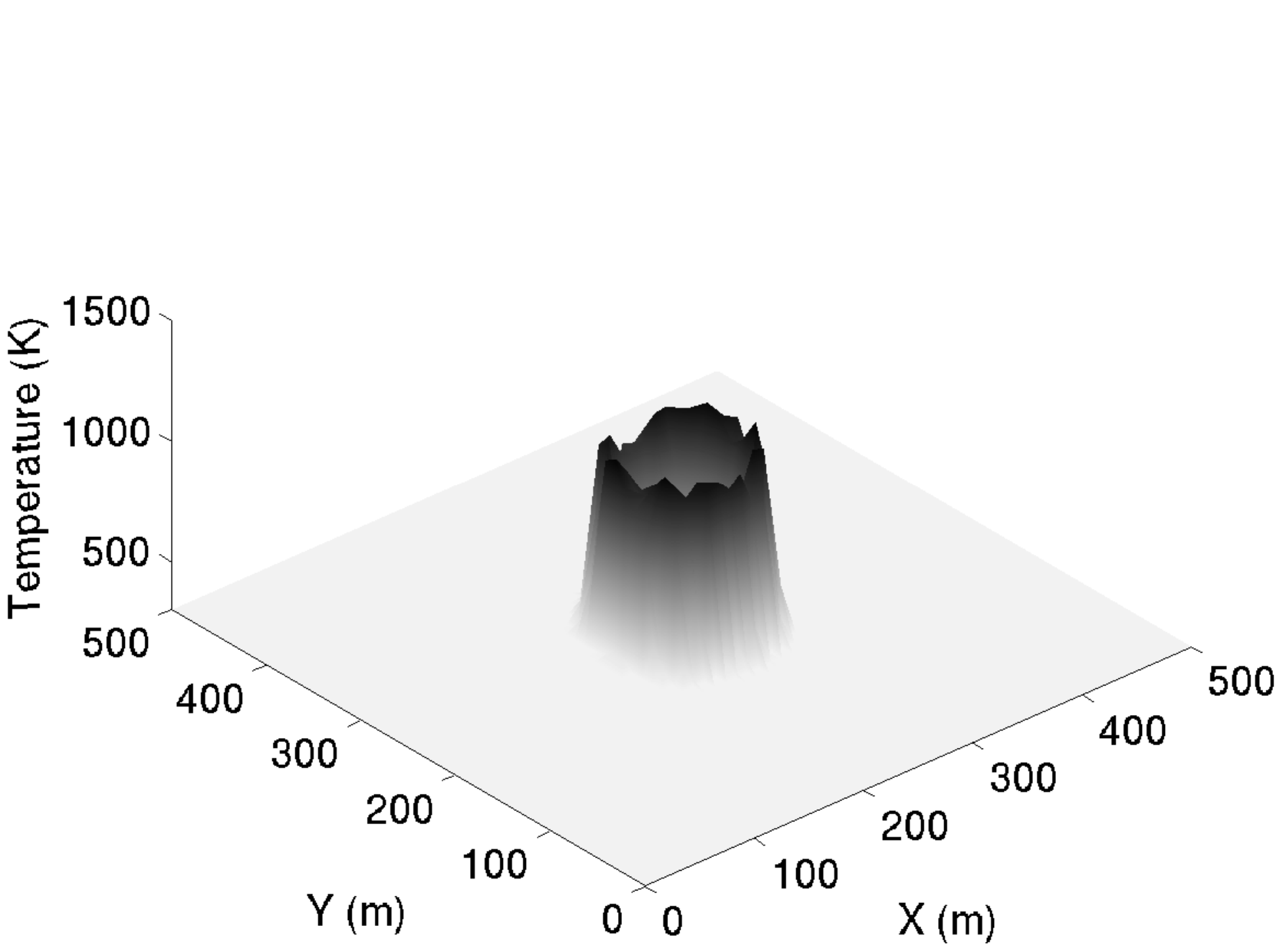}\hspace*{-1em} &
\hspace*{-1em}\includegraphics[width=1.4in]{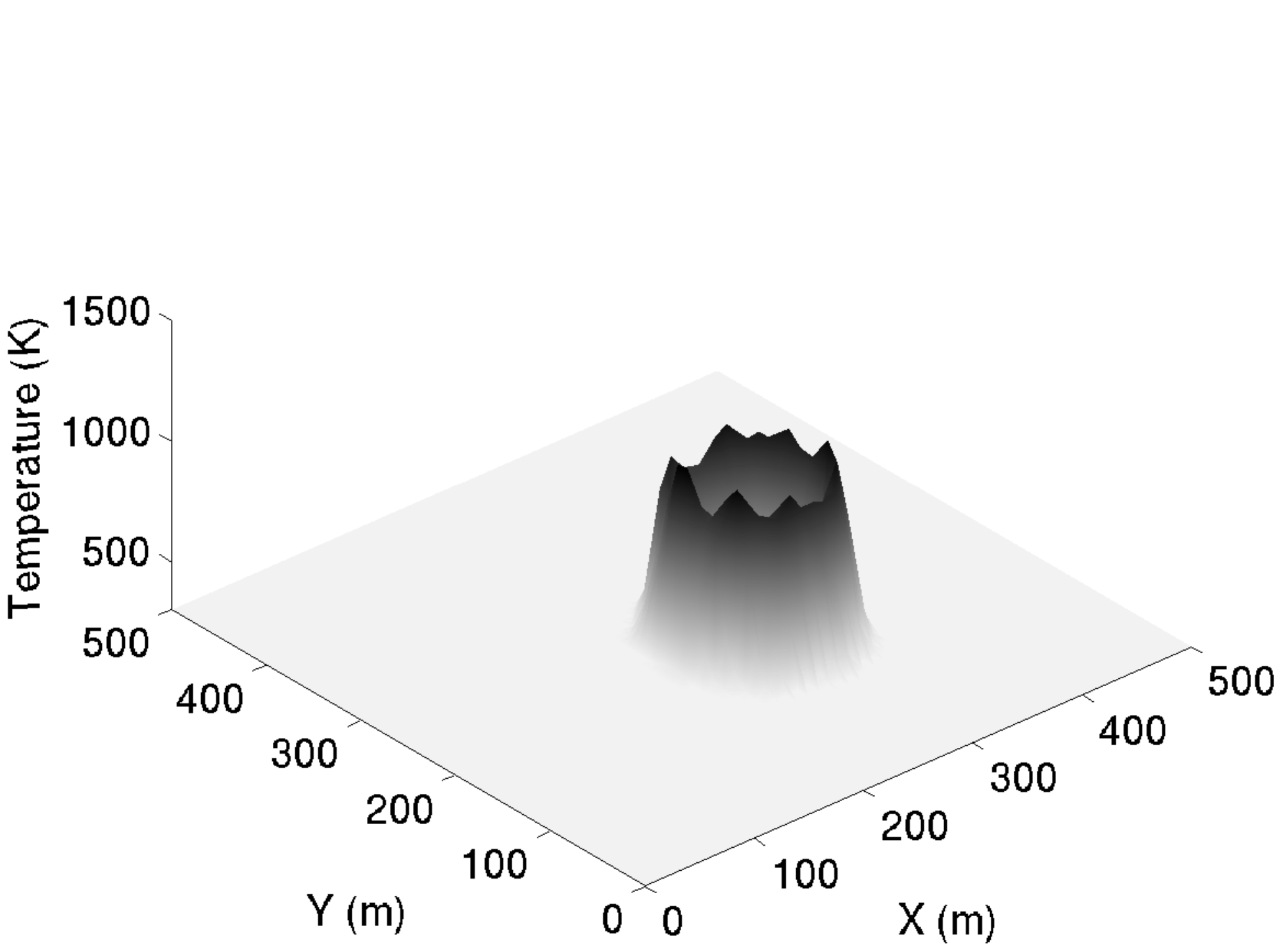}\hspace*{-1em} &
\hspace*{-1em}\includegraphics[width=1.4in]{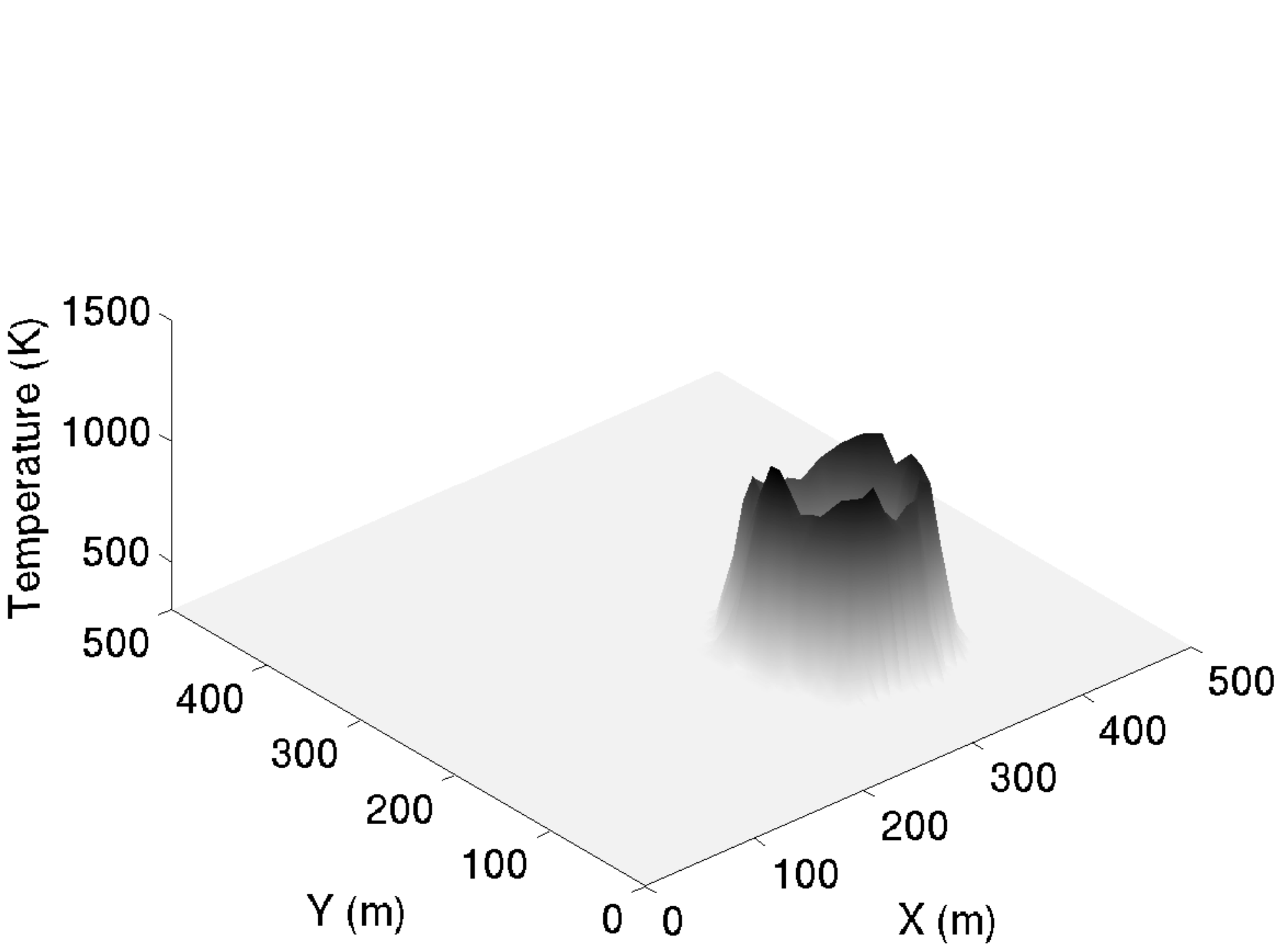}\hspace*{-1em} &
\hspace*{-1em}\includegraphics[width=1.4in]{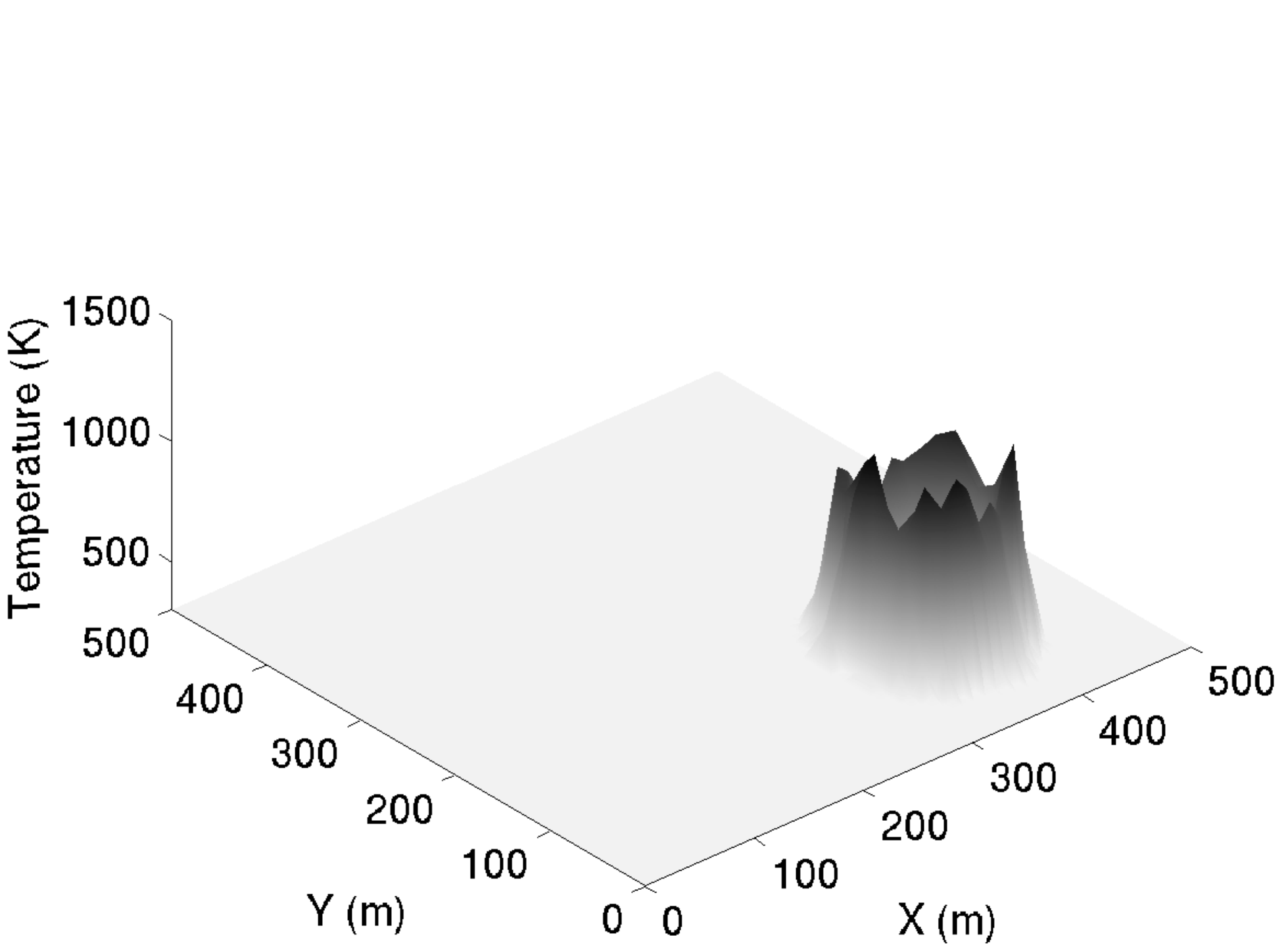}\hspace*{-1em} &
\hspace*{-1em}\includegraphics[width=1.4in]{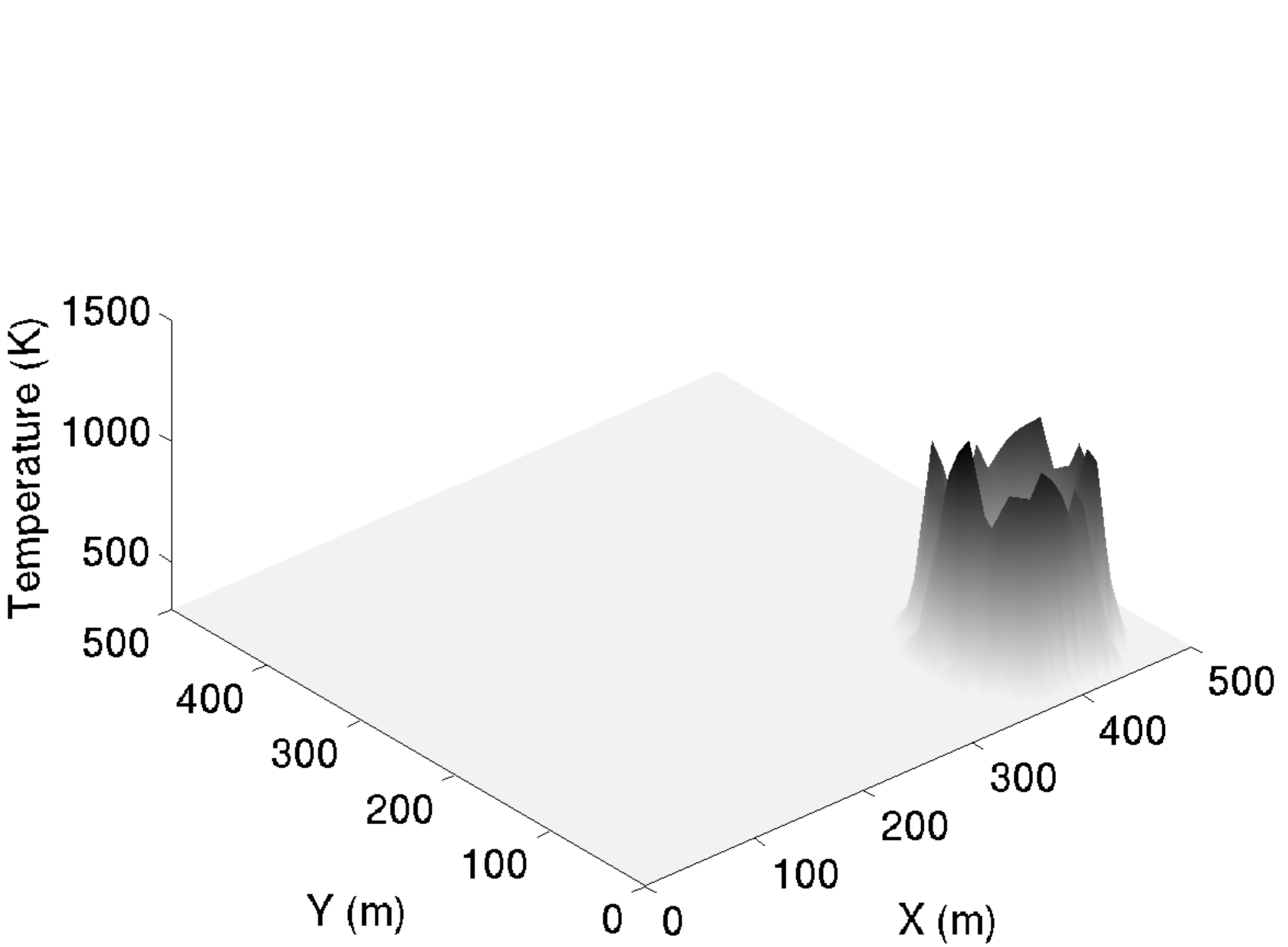}\hspace*{-1em}\\
$\lambda=0$ & $\lambda=0.25$ & $\lambda=0.5$ & $\lambda=0.75$ & $\lambda=1$%
\end{tabular}
\end{center}
\caption{\captionmorph}%
\label{fig:morph}%
\end{figure}

\newcommand{\captionsearch}{Update of the value of $I+T$ at one point
by local optimization. The mapping $I+T$ transforms the rectangle
$BCDEFGJ$ into the deformed shape $B'C'D'E'F'G'J'$. The value
$A'=(I+T)(A)$ is to be determined by minimizing the objective
function $J_{D}(T)$ defined by (\ref{eq:registration}), varying only
the location of the point $A^{\prime}$, starting from a given initial
value $A_{0}^{\prime}$. First, the objective function is evaluated at
several search points (empty circles) and then, starting from the
best location found, it is locally optimized by coordinate descent
alternating in the x and y directions, using a standard library
function; we have used \texttt{fminbnd} from MATLAB. The mapping
$I+T$ is defined by bilinear interpolation in each quadrant (e.g.,
$A^{\prime}B^{\prime}C^{\prime}D^{\prime}$), and the inverse
$(I+T)^{-1}$ on the quadrant exists if and only if the quadrant is
convex \cite{Frey-1978-SRG}. Thus, $A^{\prime}$ is constrained within
the quadrilateral $B^{\prime}D^{\prime}F^{\prime}H^{\prime }$ (thick
dot lines) formed by the midpoints of the sides the deformed
rectangle. The search points are constructed by putting several local
grid points (here, $2$) at equal distance on the segment between the
initial position $A_{0}^{\prime}$ and the midpoint (e.g.,
$A_{0}^{\prime}B^{\prime}$), and then adding equidistant grid of
search points on each of the lines that connect corresponding points,
to form a triangle as shown.}
\begin{figure}
\begin{center}
\includegraphics[width=4in]{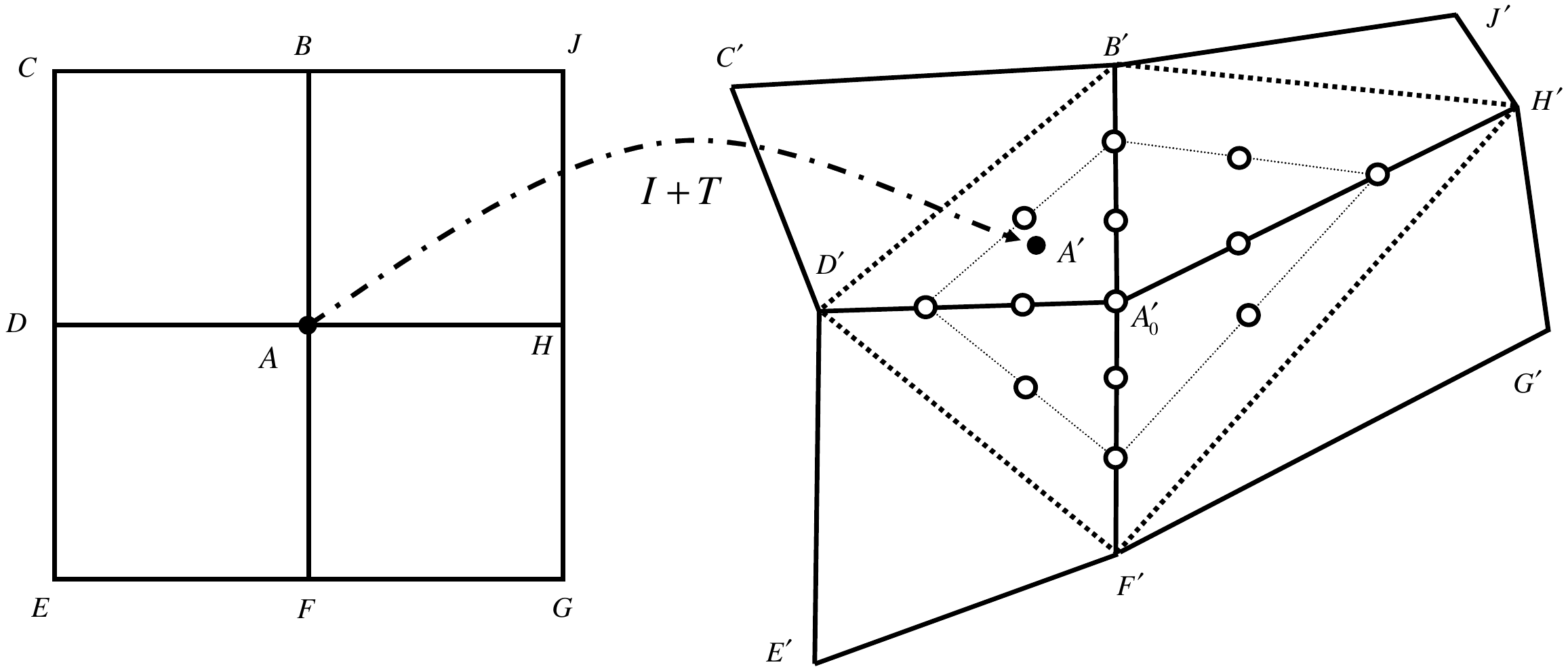}
\end{center}
\caption{\captionsearch}%
\label{fig:search}%
\end{figure}

\newcommand{\captiondapdf}{Probability densities estimated by a~Gaussian kernel with bandwidths
$37$ $K$, $19$ $K$, and $30$ $m$. Data was collected from the
ensemble shown in Fig.~\ref{fig:da_morph}c. Displayed are typical
marginal probability densities at a point near the reaction area of the original
temperature (a), the registration residual of temperature (b), and the
registration mapping component in the $x$ direction (c). The
transformation has vastly improved Gaussian nature of the densities.}
\begin{figure}
\begin{center}%
\begin{tabular}
[c]{ccc}%
\includegraphics[width=1.5in]{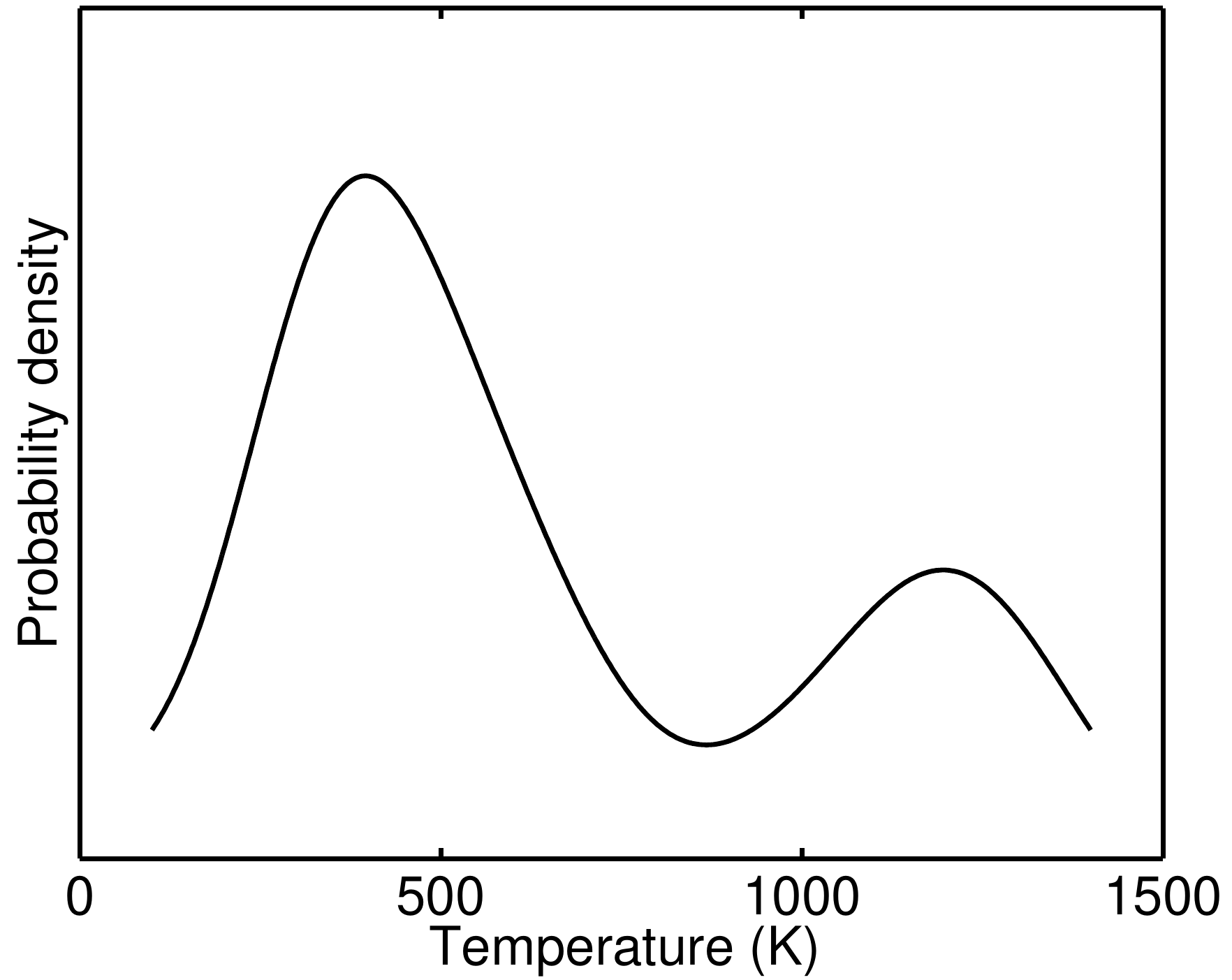} &
\includegraphics[width=1.5in]{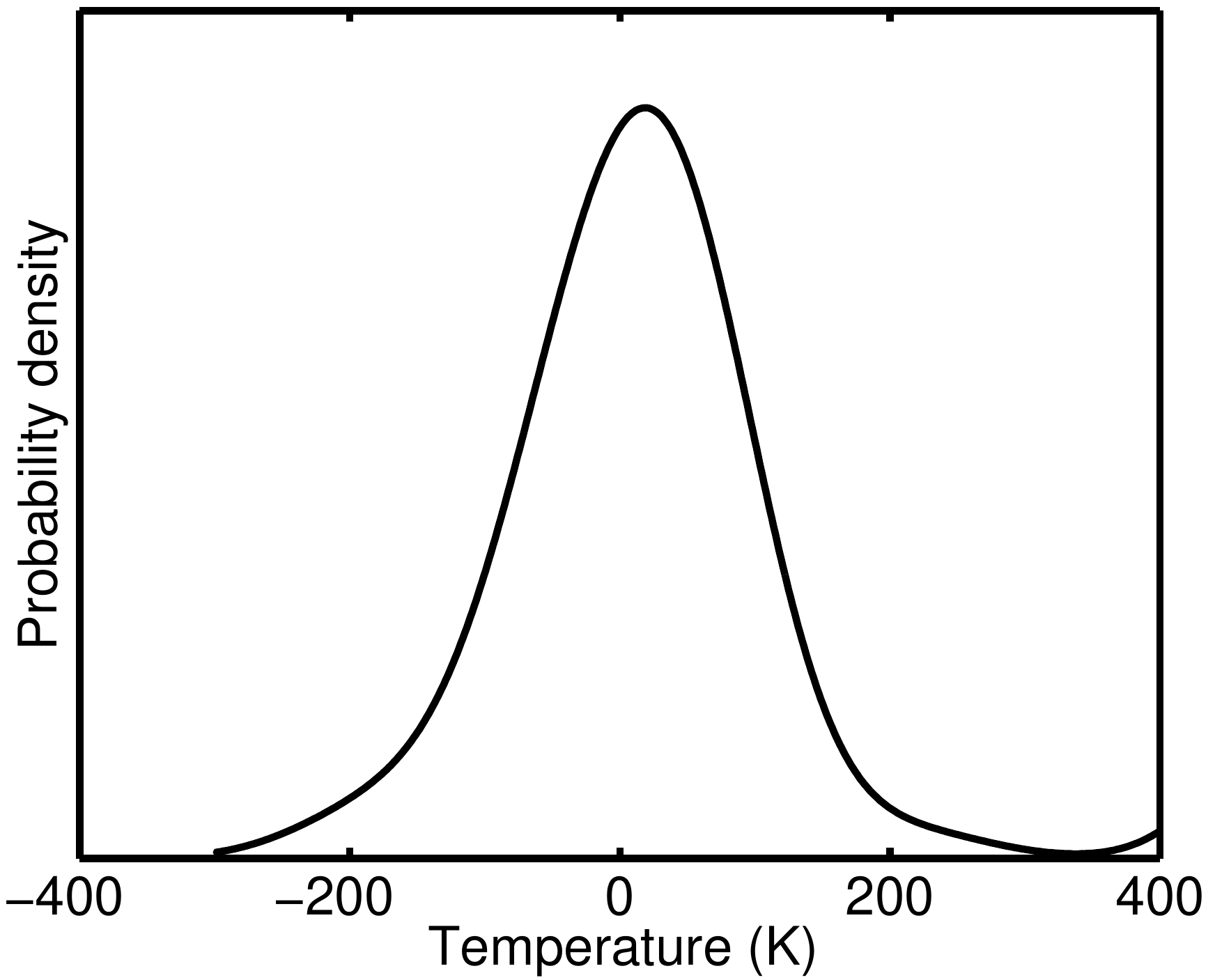} &
\includegraphics[width=1.5in]{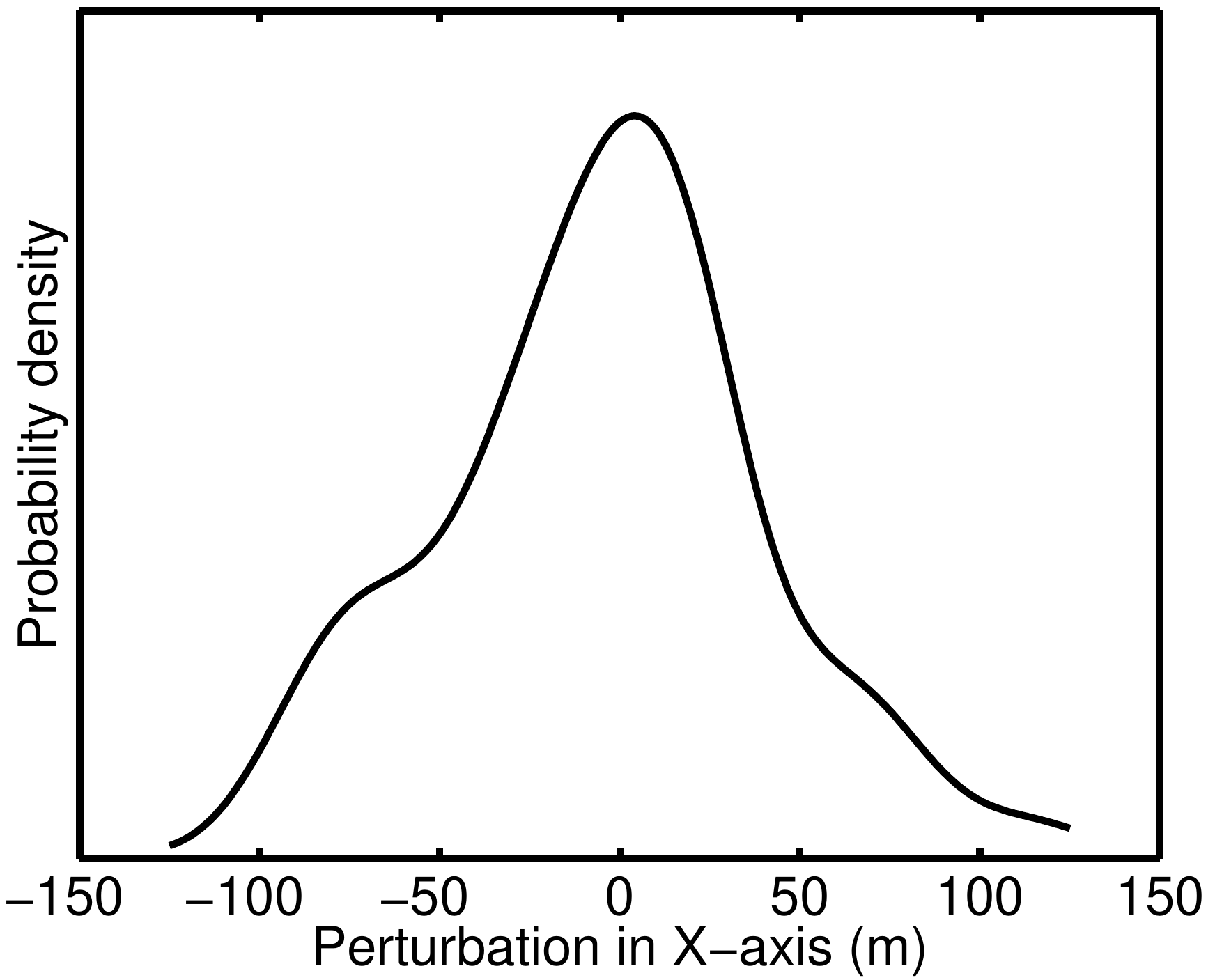}\\
(a) & (b) & (c)
\end{tabular}
\end{center}
\caption{\captiondapdf}%
\label{fig:da_pdf}%
\end{figure}

\newcommand{\captiondamorph}{Data assimilation by the morphing ensemble filter. The forecast
ensemble (b) was created by smooth random morphing of the initial
temperature profile (a). The analysis ensemble (d) was obtained by
the EnKF applied to the transformed state. The data for the EnKF was
the morphing transformation of the simulated data (c), and the
observation function was the identity mapping. Contours are at $800
\,K$, indicating the location of the fireline. The reaction zone is
approximately between the two curves.}
\begin{figure}
\begin{center}%
\begin{tabular}
[c]{cccc}%
\includegraphics[width=1.4in]{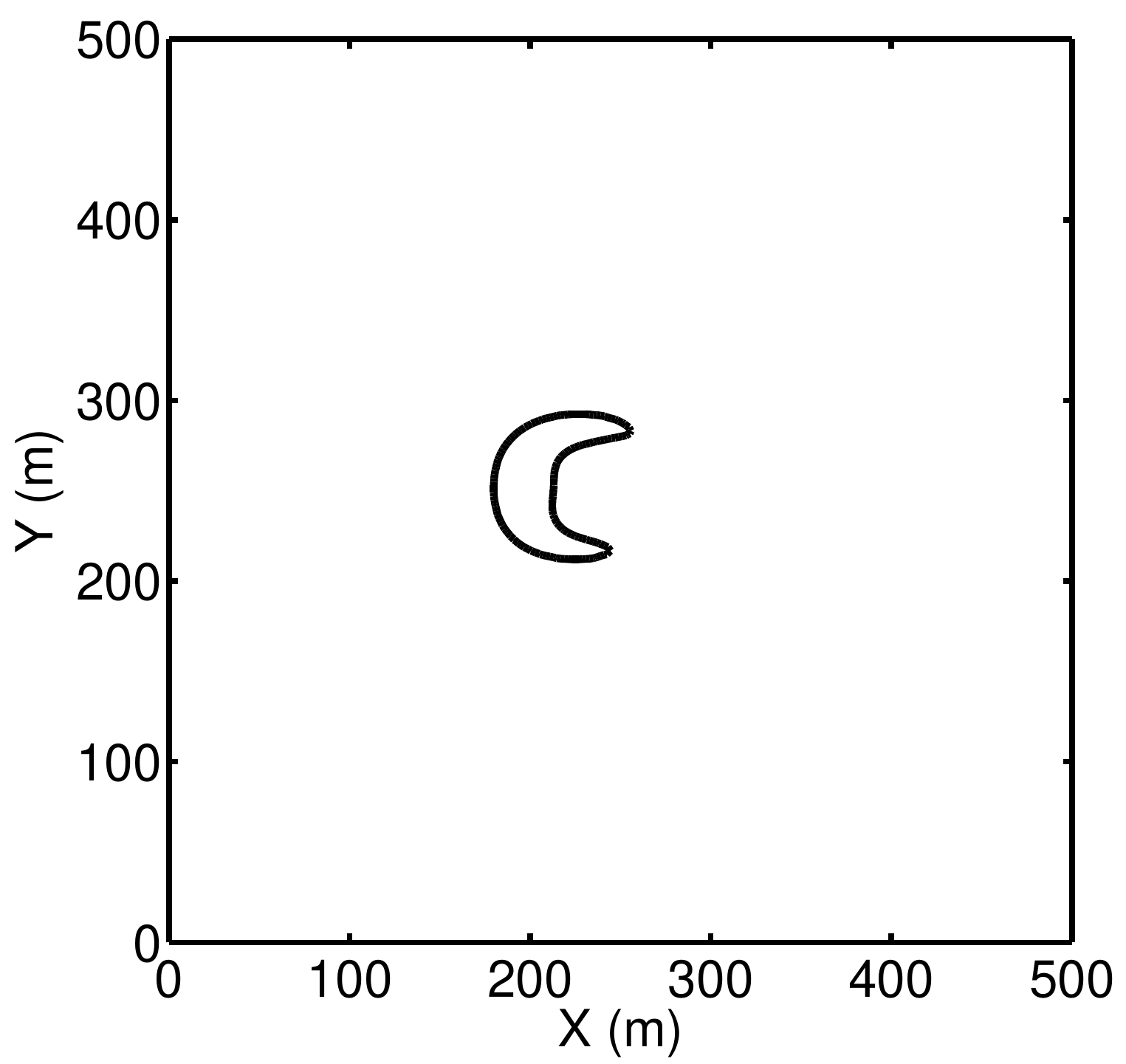} &
\includegraphics[width=1.4in]{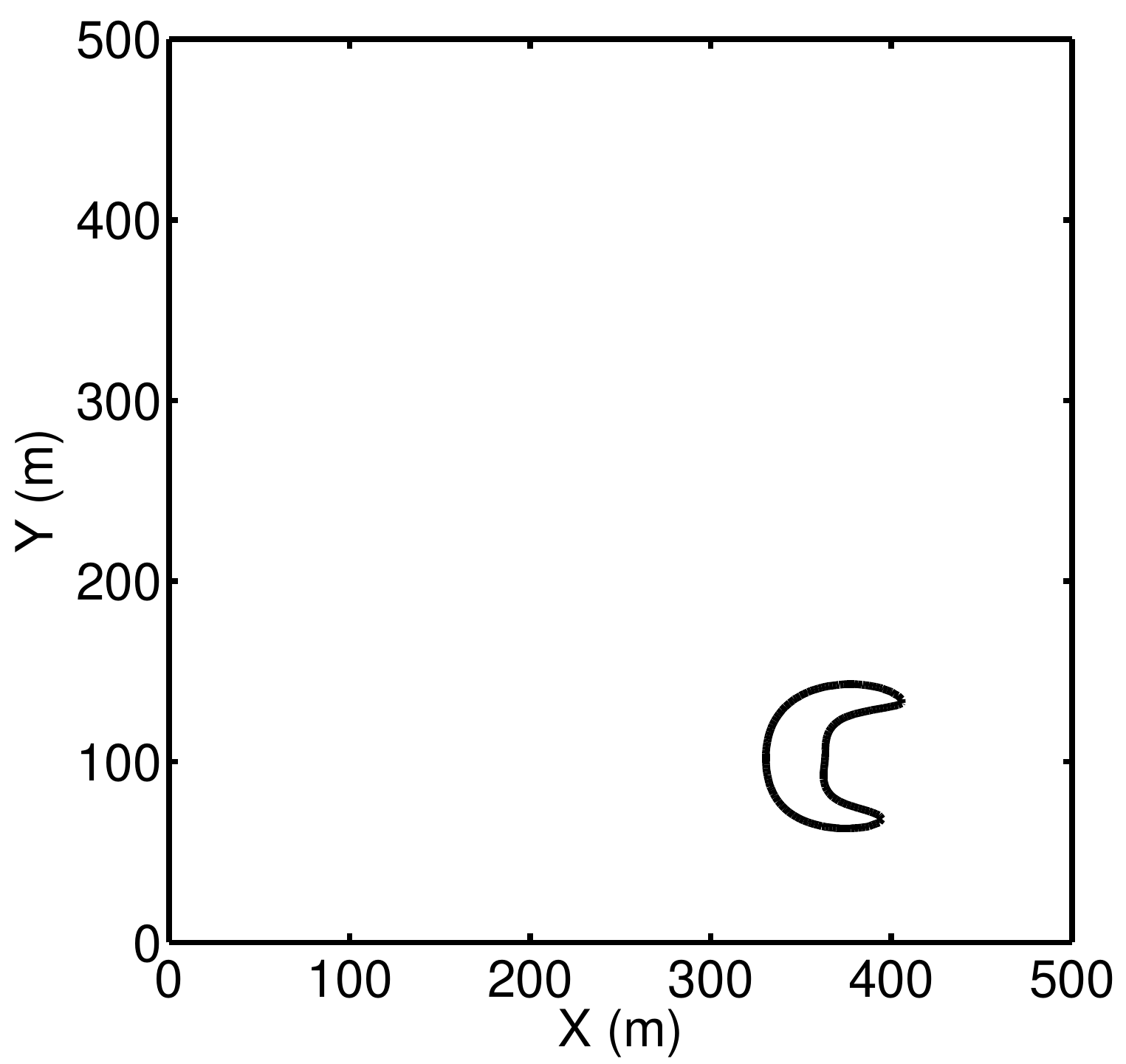} &
\includegraphics[width=1.4in]{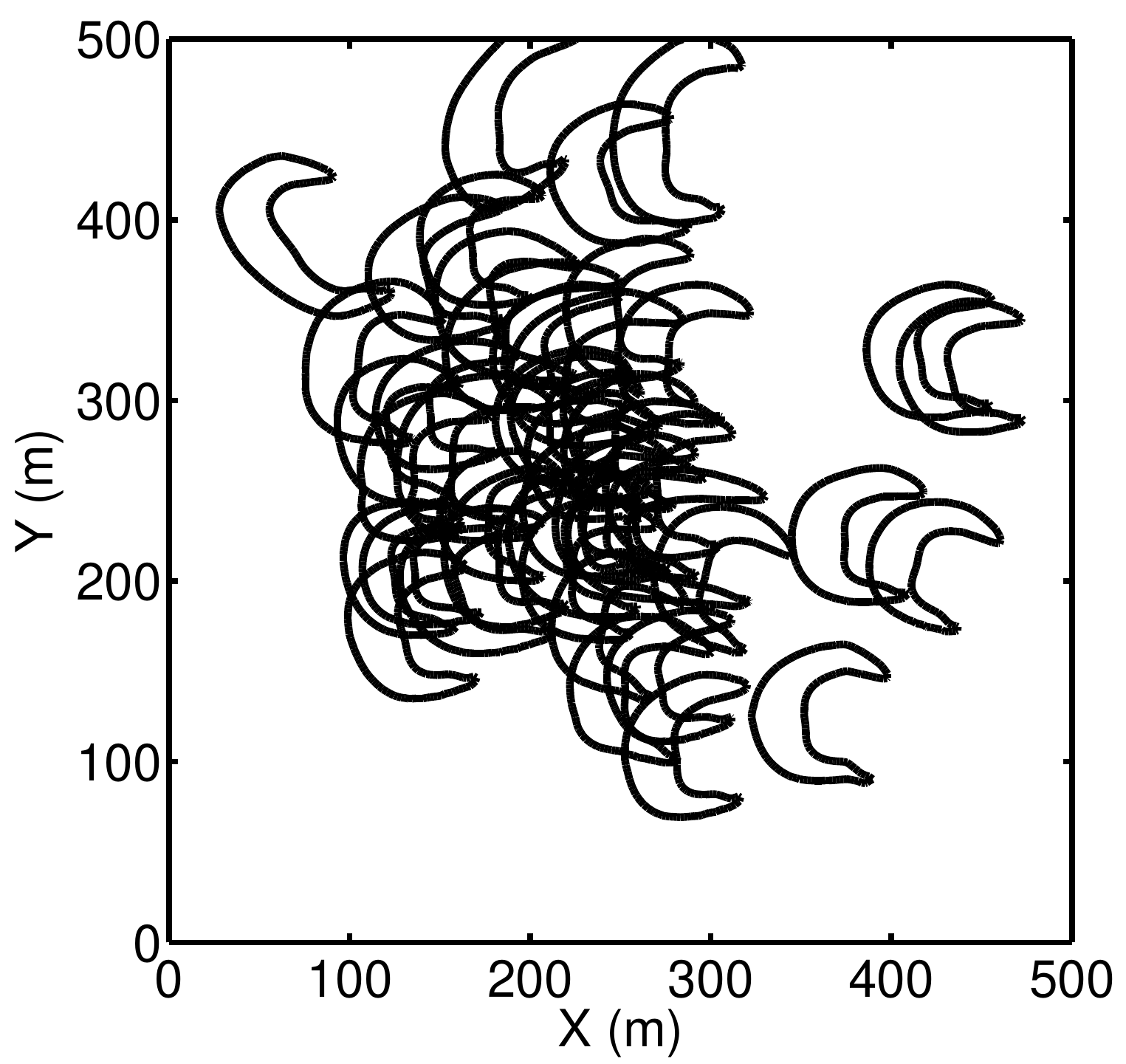} &
\includegraphics[width=1.4in]{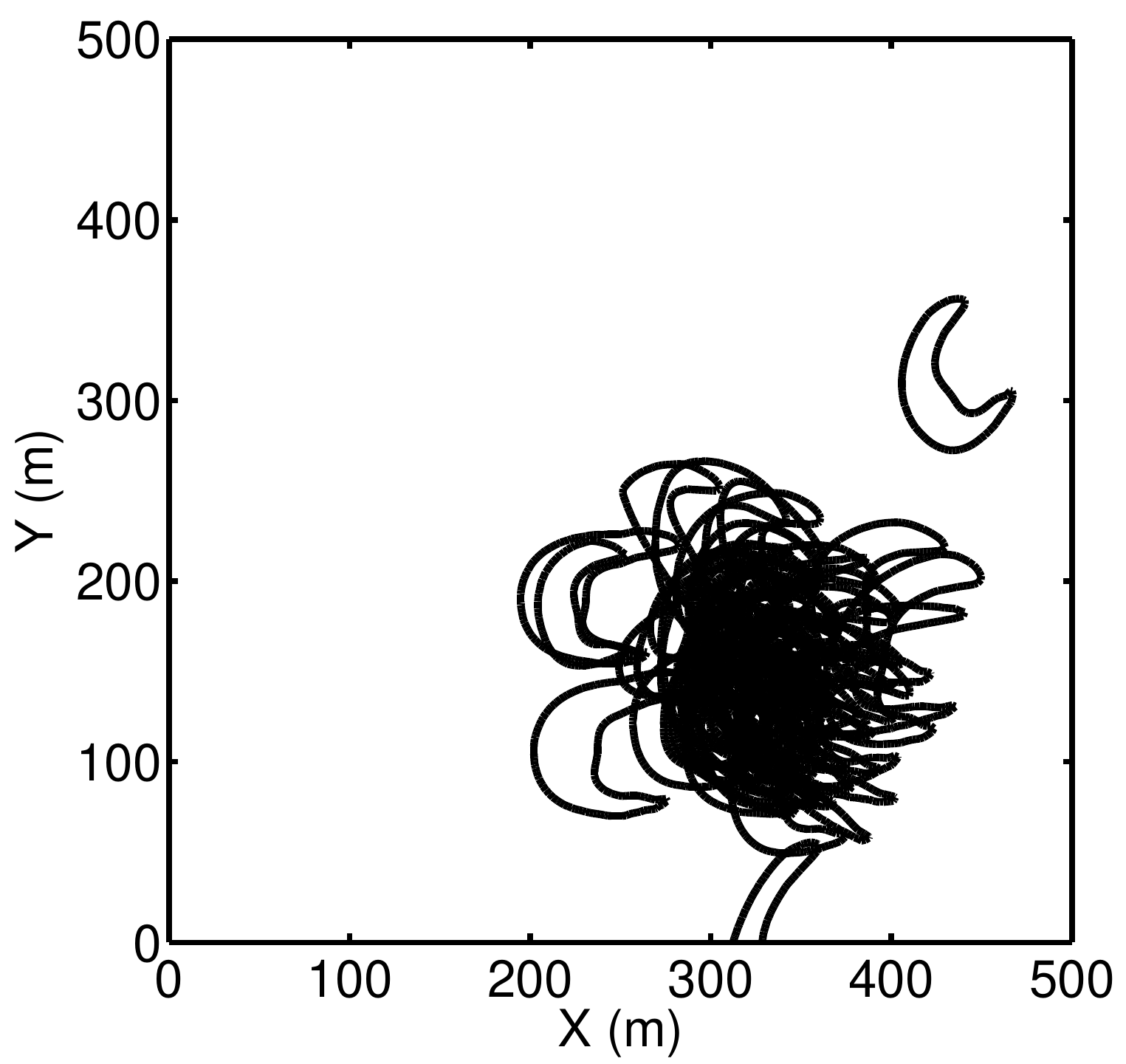}\\
(a) & (b) & (c) & (d)
\end{tabular}
\end{center}
\caption{\captiondamorph}%
\label{fig:da_morph}%
\end{figure}

\newcommand{\captiondamorphfive}{After five analysis cycles, the ensemble shows less spread and
follows the data reliably. Ensemble members were registered using the
initial state, advanced in time without data assimilation (a). The forecast ensemble (b) is closer to the
simulated data (c) because of preceding analysis steps that have
attracted the ensemble to the truth. The analysis ensemble (d) has a
little less spread than the forecast, and the change between the
forecast and the analysis is well within the ensemble spread.
Contours are at $800\, K$, indicating the location of the fireline. The
reaction zone is approximately between the two curves.}
\begin{figure}
\begin{center}%
\begin{tabular}
[c]{cccc}%
\includegraphics[width=1.4in]{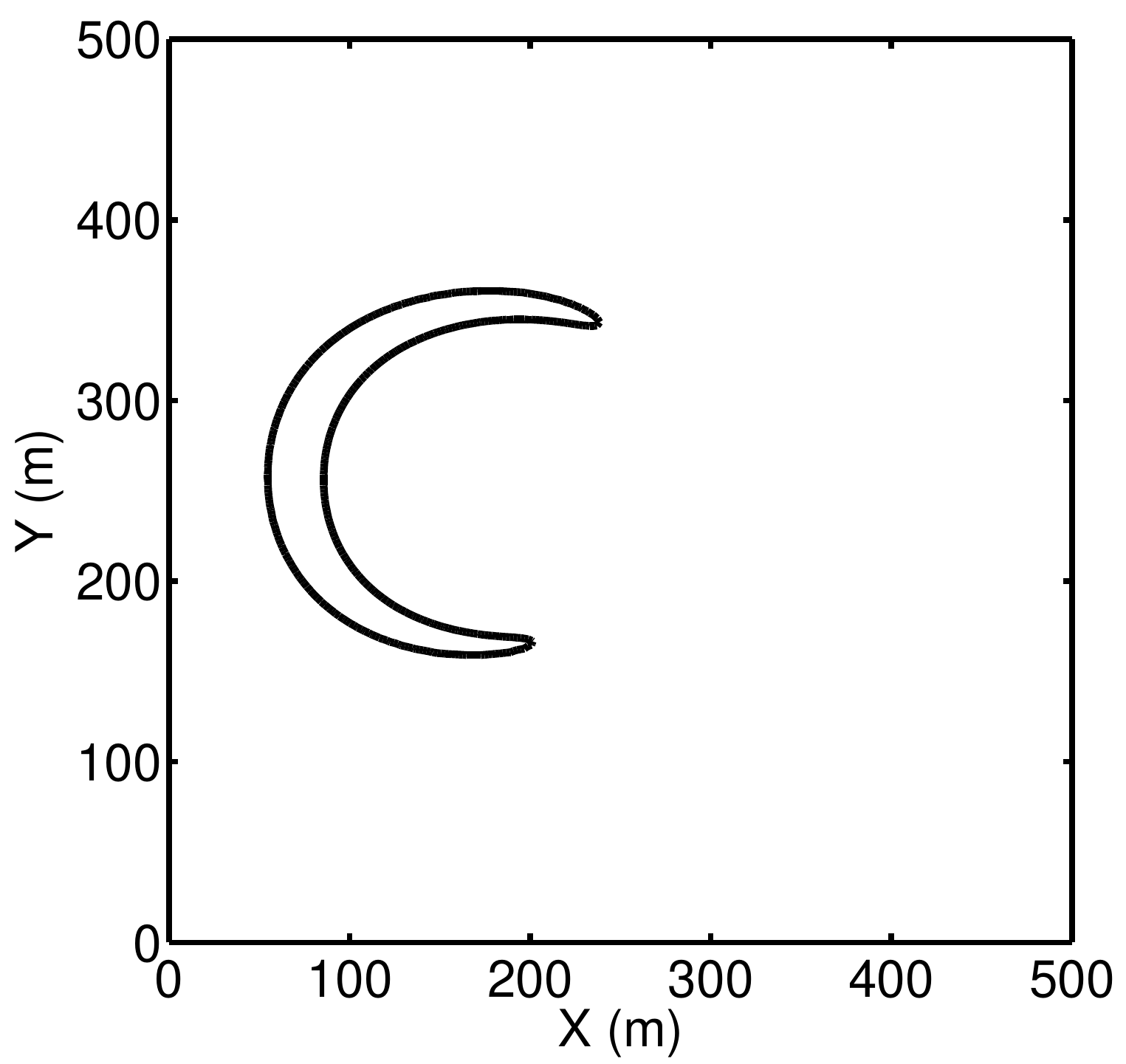} &
\includegraphics[width=1.4in]{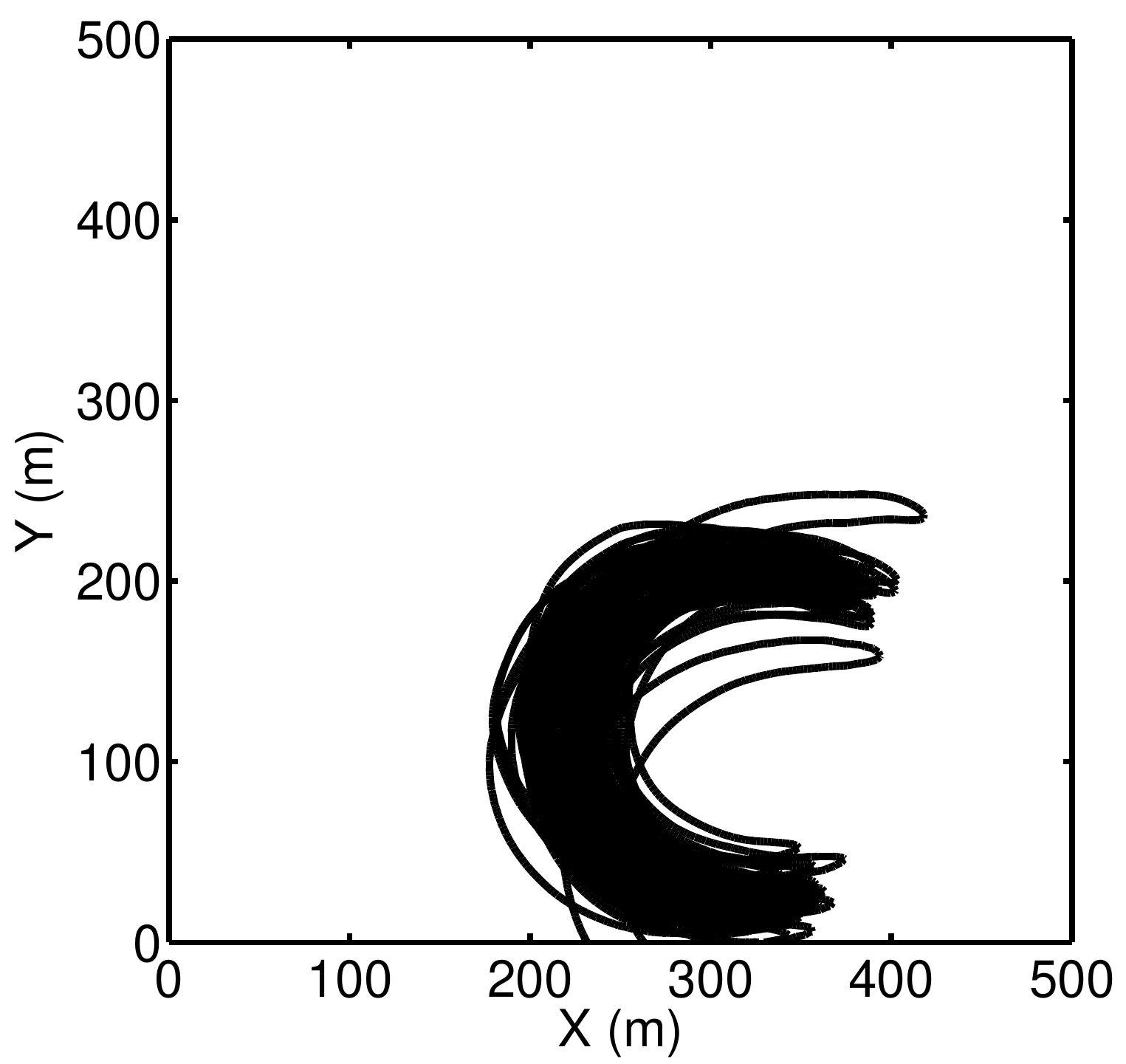} &
\includegraphics[width=1.4in]{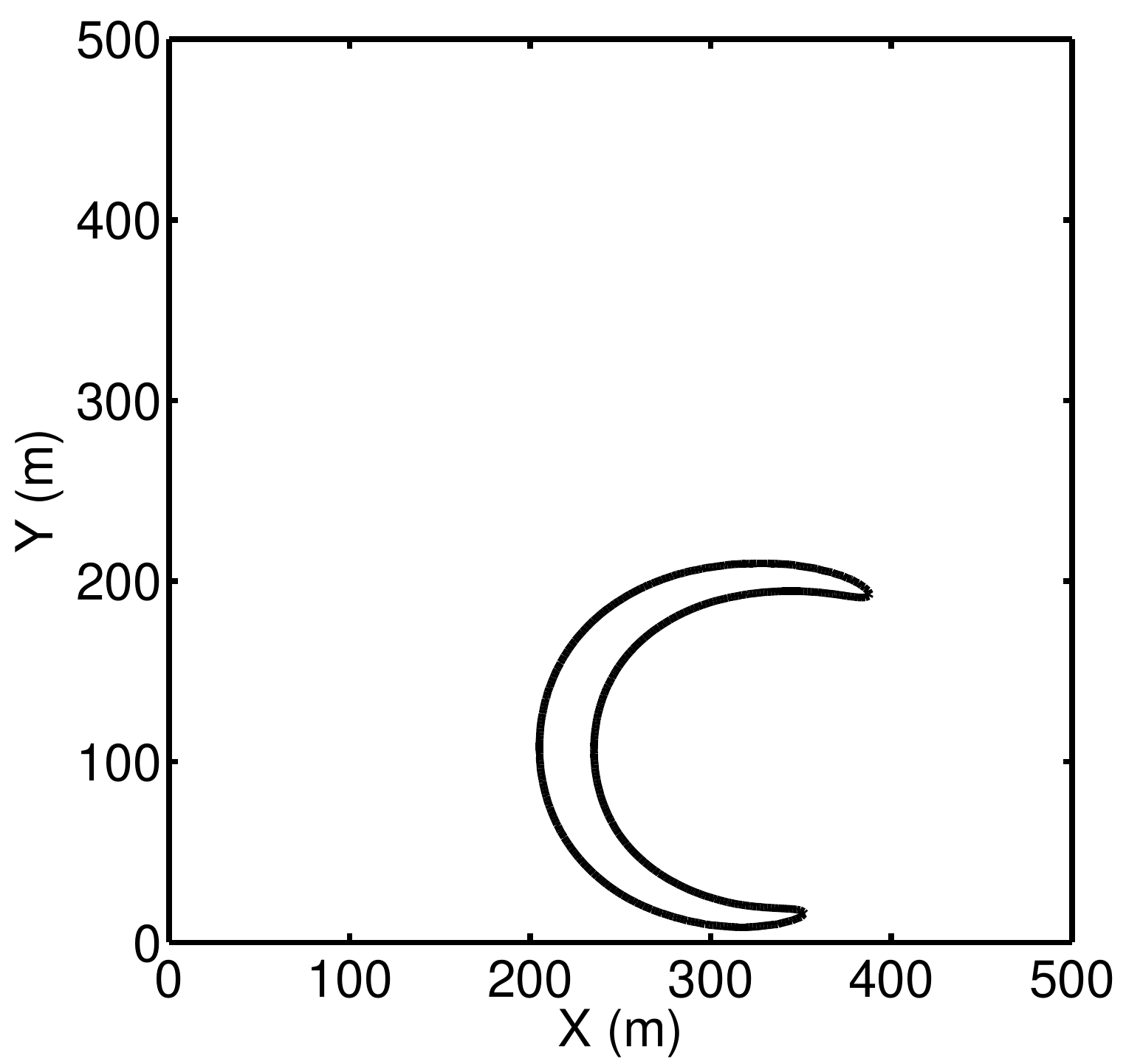} &
\includegraphics[width=1.4in]{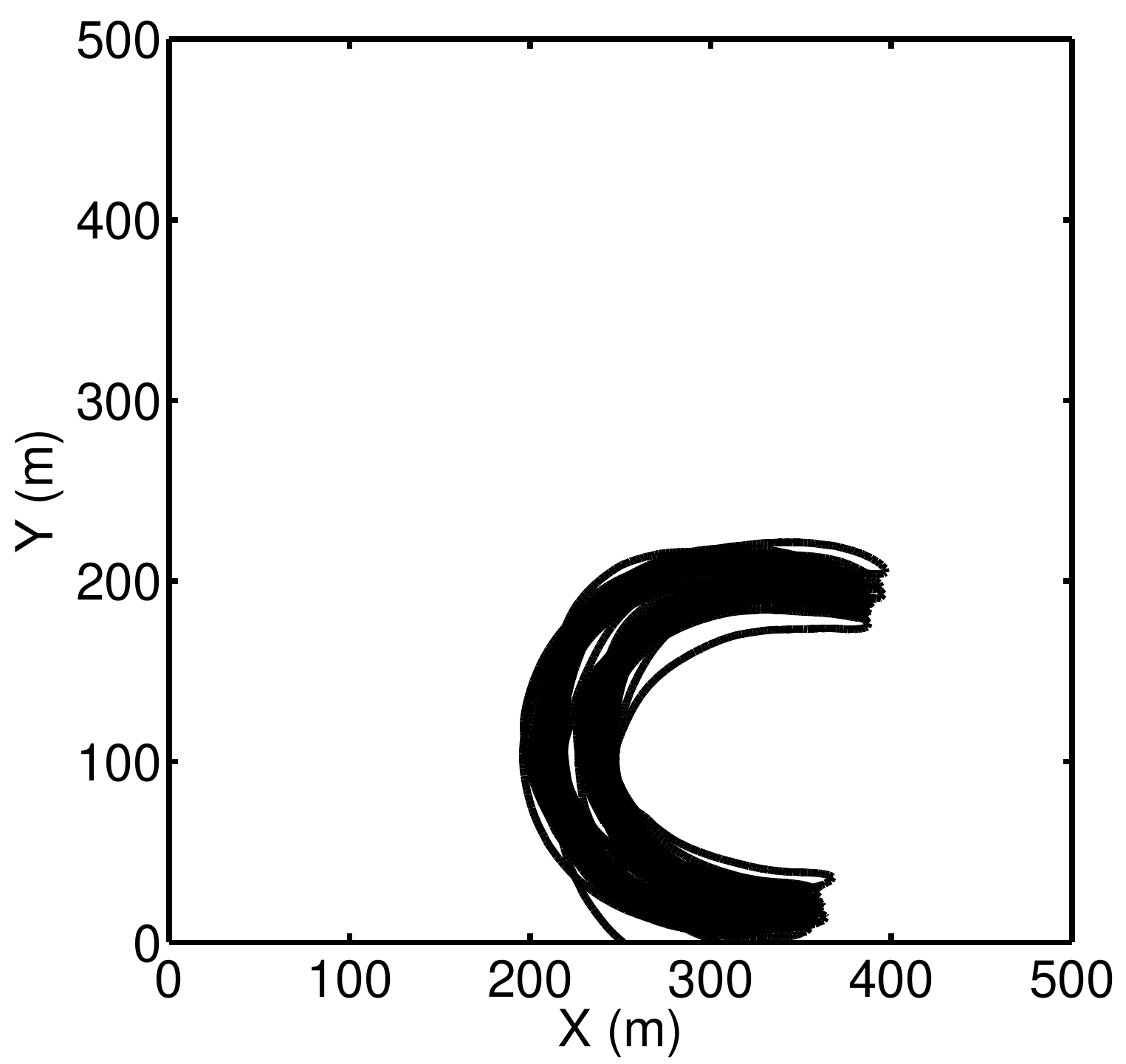}\\
(a) & (b) & (c) & (d)
\end{tabular}
\end{center}
\caption{\captiondamorphfive}%
\label{fig:da_morph5}%
\end{figure}

\newcommand{\captiondanormtestfive}{
The $p-$value from the Anderson-Darling test of the data from the
ensemble after five morphing EnKF analysis cycles shows the
ensemble transformed into its registration representation, the registration
residual of the temperature (b) and the registration mapping (c), has
distribution much closer to Gaussian than the original ensemble (a).
The shading at a point indicates where the marginal distribution of the ensemble at
that point is highly Gaussian (white) or highly non-Gaussian (black).}
\begin{figure}
\begin{center}%
\begin{tabular}
[c]{ccc}%
\includegraphics[width=1.5in]{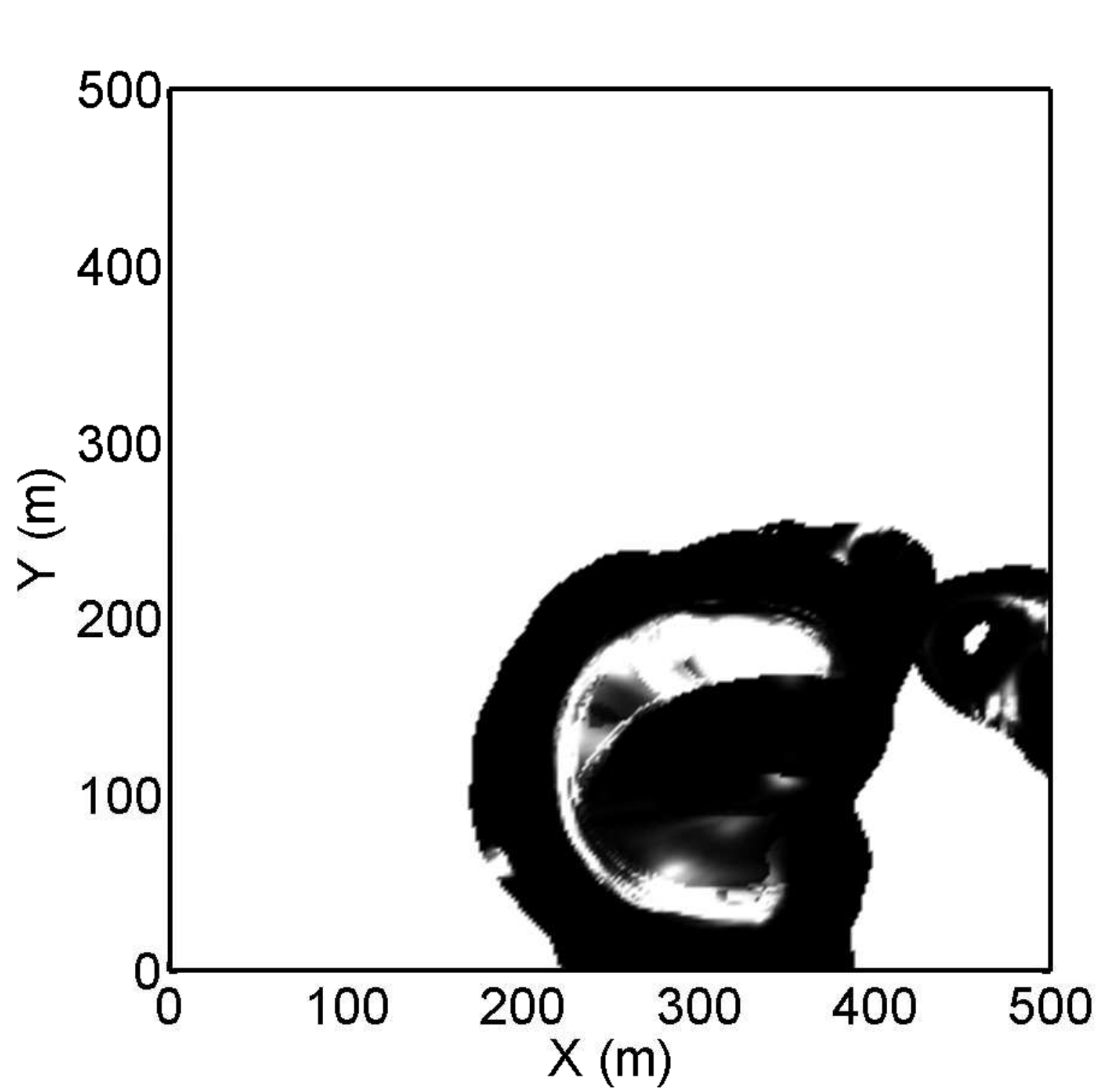} &
\includegraphics[width=1.5in]{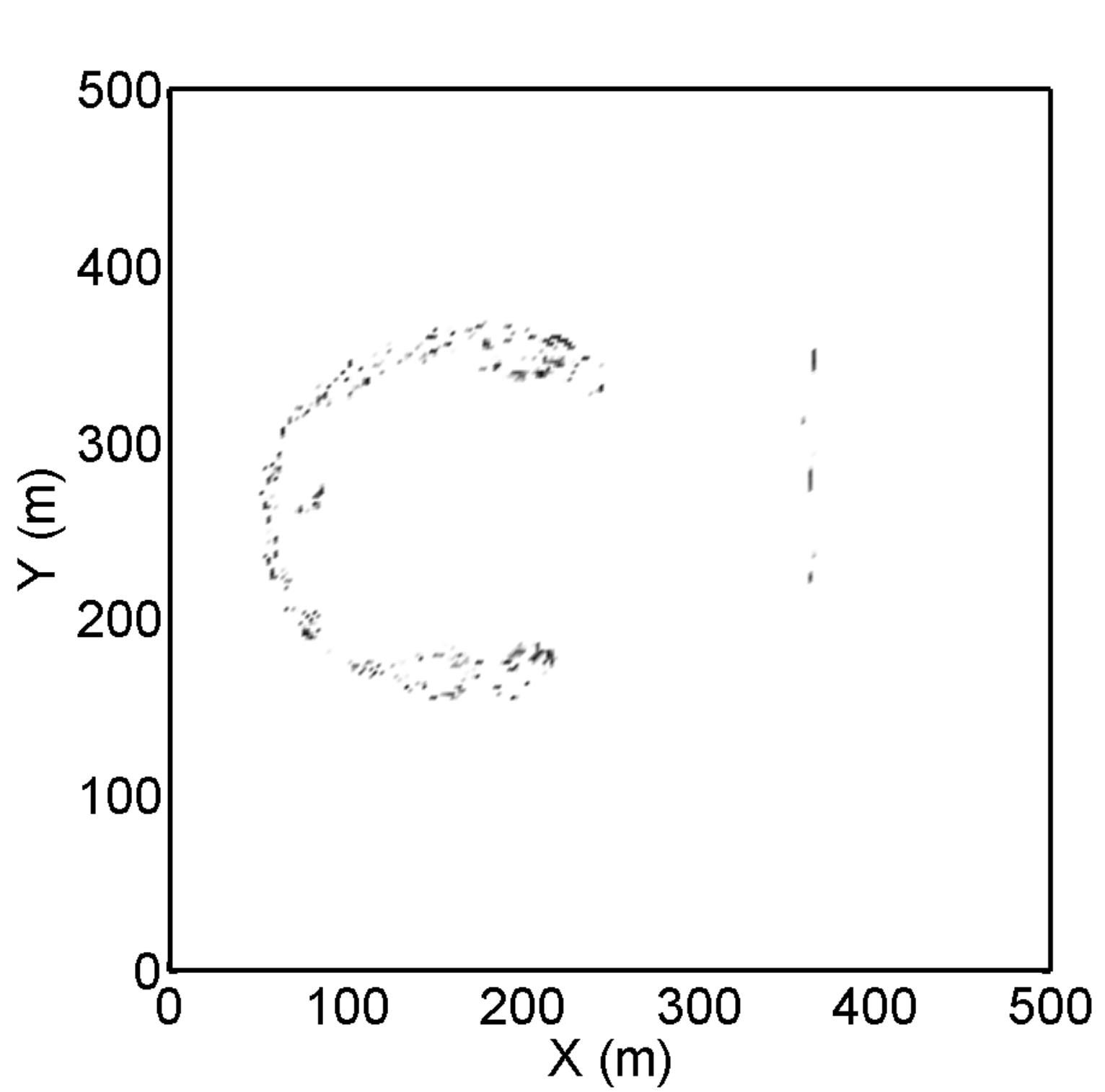} &
\includegraphics[width=1.5in]{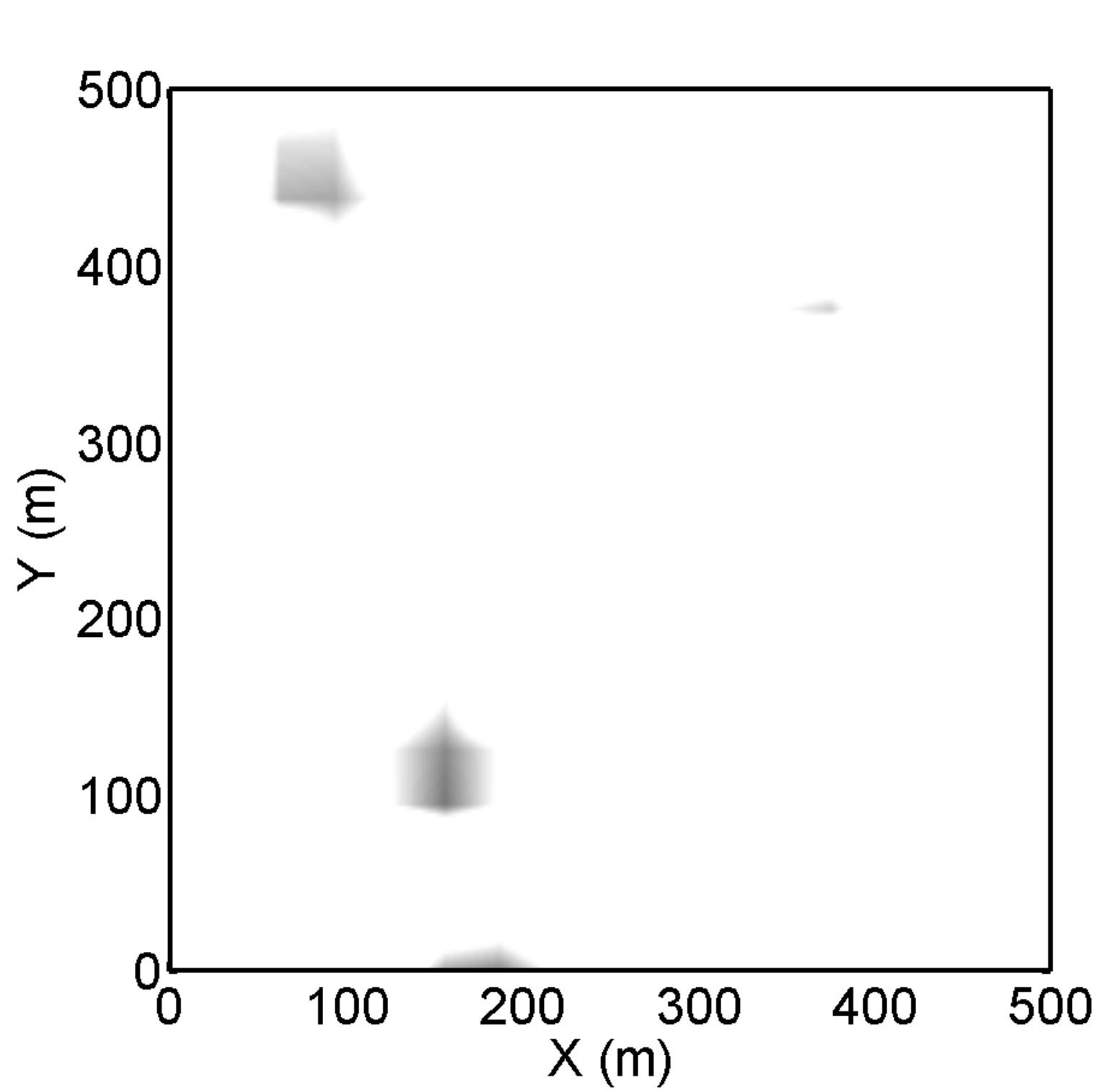}\\
(a) & (b) & (c)
\end{tabular}
\end{center}
\caption{\captiondanormtestfive}%
\label{fig:da_normtest5}%
\end{figure}

\end{document}